\documentclass{amsart}
\usepackage{amsmath}
\usepackage{amsfonts}
\usepackage{mathrsfs}
\usepackage{bbm}
\usepackage{bm}

\usepackage[normalem]{ulem}
\usepackage{enumerate}

\newcommand{\bpf}[1][Proof]{{\noindent {\sc #1: }}}
\newcommand{\epf}{{{\hfill $\Box$ \smallskip}}}

\newcommand{\R}{\mathbb{R}}
\newcommand{\N}{\mathbb{N}}

\newcommand{\tbf}{\mathbf{t}}

\newcommand{\sbf}{\mathbf{s}}

\newcommand{\ibf}{\mathbf{i}}
\newcommand{\lambdabf}{\mathbf{\lambda}}

\newcommand{\Pp}{\mathsf{P}}
\newcommand{\pp}{\mathsf{p}}

\newcommand{\K}{\mathcal{K}}
\newcommand{\Sc}{\mathcal{S}}
\newcommand{\Rc}{\mathcal{R}}
\newcommand{\Tc}{\mathcal{T}}

\newcommand{\Is}{\mathscr{I}}

\newcommand{\Cs}{\mathscr{C}}

\newtheorem{theorem}{Theorem}
\newtheorem{lemma}{Lemma}
\newtheorem{remark}{Remark}
\newtheorem{corollary}{Corollary}
\newtheorem{proposition}{Proposition}
\newtheorem{example}{Example}

\begin{document}

\title[Regularity of Invariant Densities]{Regularity of invariant
densities for 1D-systems with random switching}

\author{Yuri Bakhtin}
\address{Yuri Bakhtin \\
School of Mathematics \\
Georgia Institute of Technology\\
Atlanta GA \\ 
USA}

\author{Tobias Hurth}
\address{Tobias Hurth \\
School of Mathematics \\
Georgia Institute of Technology\\
Atlanta GA \\ 
USA}

\author{Jonathan C. Mattingly}
\address{Jonathan C. Mattingly \\
Mathematics Department \\
Center of Theoretical and Mathematical Science \\
Department of Statistical Science \\
and Center of Nonlinear and Complex Systems \\
Duke University \\
Durham NC \\
USA}

\begin{abstract}
This is a detailed analysis of invariant measures for one-dimensional dynamical
systems with random switching. In particular, we prove smoothness of
the invariant densities away from critical points and describe the asymptotics
of the invariant densities at critical points.
\end{abstract}

\maketitle

\bigskip



\section{Introduction}

In this paper, we study the regularity theory for invariant densities of
dynamical systems with random switching (switching systems, in short) with
one-dimensional continuous
component. Dynamical systems with random switching are also
known as piecewise deterministic Markov processes.  

We show that smoothness of the vector fields
governing the
dynamics translates into smoothness of invariant densities away from
critical points of the vector fields (Theorem~\ref{thm:smoothness}). For
analytic vector fields, we derive the asymptotically dominant term of an
invariant density as its argument approaches a critical point of the
corresponding vector field (Theorem~\ref{thm:critical_points_analytic}).  

In the literature, regularity properties of
invariant densities are often assumed in order to derive other features of the
densities. For instance, it is shown in~\cite[Proposition 3.1]{Faggionato} that
if invariant densities are $\Cs^1$
on a set $\Omega$, they satisfy the Fokker--Planck equations associated to
the switching system in the interior of $\Omega$. From this differential
characterization,
the authors deduce time-reversibility of stationary piecewise deterministic
Markov processes and derive explicit formulas
for the invariant densities of certain switching systems that they call exactly
solvable. A result similar to~\cite[Proposition 3.1]{Faggionato} can be found
in~\cite[Theorem 1]{Horton}. In this paper, we give sufficient conditions for
continuity and differentiability of invariant densities that are stated in
terms of the vector fields, and are easily verifiable. In particular, we show
that if none of the vector fields vanish at a point $\xi \in
\R$ and if all vector fields are $\Cs^{n+1}$ in a neighborhood of $\xi$, then
the invariant densities are $\Cs^{n}$ at $\xi$. This statement is not
surprising, but as far as we know, no rigorous proof of it has been given yet. 

The second question addressed in this paper is how invariant densities behave
at critical points of the vector fields. In the case of two vector fields on a
bounded interval that point in opposite directions, \cite[Proposition
3.12]{Faggionato} gives
an explicit formula for the invariant densities. From this formula, one
obtains the exact asymptotic behavior of the densities close to
critical points. However, 
computing invariant densities explicitly is in general very difficult
(\cite[Section 3.3]{Faggionato}). Finding necessary and sufficient conditions
for boundedness of invariant densities is already challenging. In the
one-dimensional
case, invariant densities are bounded away from critical
points (Lemma~\ref{lm:bounded}), but we expect to find switching systems with
two-dimensional continuous component whose invariant
densities become
unbounded along curves that do not contain any critical points.
If the continuous component is one-dimensional, \cite[Theorem 1]{Balazs}
provides sufficient conditions for boundedness of an invariant density close to
a critical point of its associated vector field. For vector fields that behave
linearly close to a critical point, we give necessary conditions and sufficient
conditions for boundedness
in terms of the vector fields and the switching rates
(Corollary~\ref{co:boundedness}). These conditions recover part of the
results in~\cite{Balazs}. For analytic vector fields, we also compute the
asymptotically dominant term of an invariant density as its argument approaches
a critical point of the associated vector field
(Theorem~\ref{thm:critical_points_analytic}). Even if the vector fields are
not analytic, we can derive the asymptotically dominant term in certain cases
(Theorem~\ref{thm:critical_points}).
  
The basic tools in our investigation are two integral equations satisfied by
invariant densities. These equations are closely related to the
Fokker--Planck equations (see Appendix~\ref{sec:Fokker_Planck}), but do not
require
differentiablity of the
densities.
When deriving the asymptotically dominant
terms in the case of analytic vector fields, we use the theory of
regular singular points for systems of linear ordinary differential equations.
We follow~\cite[Section 3.11]{Taylor_ode}.  

The paper is organized as follows.
In Section~\ref{sec:definitions}, we describe the class of switching systems we
consider and
introduce some notation and terminology needed in the rest of the paper. Our
main results are stated in Section~\ref{sec:results}. In
Section~\ref{sec:representation}, we formulate the integral equations mentioned
earlier. One of these equations plays an important role in
the proof of
Theorem~\ref{thm:smoothness} that can be found in Section~\ref{sec:smoothness}. 
In Section~\ref{sec:support}, we describe the
support of invariant measures for switching systems with one-dimensional
continuous
component. We exploit this description of the support in
Section~\ref{sec:critical_points}, which is devoted to proving
Theorems~\ref{thm:critical_points_analytic}
and~\ref{thm:critical_points}. Section~\ref{sec:integral_equations} contains
the proofs of the integral equations from Section~\ref{sec:representation}.
The appendix contains a remark on how the integral equations relate to the
Fokker--Planck equations for the invariant densities.  

\subsection*{Acknowledgements} JCM acknowledges the support of the
National Science Foundation(DMS-0854879) and the Simons Foundation.

\section{Definitions and Notation}  \label{sec:definitions}
We consider a similar setup as in~\cite{Bakhtin}, but we restrict our analysis
to switching systems with one-dimensional continuous component.  
Let $D$ be a finite collection of continuously differentiable and 
forward-complete vector fields on
$\R$. We denote these vector fields by $u_i$, where $i$ is some element of
the index set $S = \{1, \ldots, n\}$. 

Our assumptions on the vector fields imply
that the initial-value problem 
\begin{align*}
 \dot{x}(t) &= u_i(x(t)), \\
x(0)        &= \xi
\end{align*}
induced by $u_i \in D$ and by $\xi \in \R$ has a unique solution defined for all
$t \geq 0$. We define a stochastic process $X = (X_t)_{t \geq 0}$ on $\R$ as
follows.
Given an initial point $\xi \in \R$ and an initial driving vector field $u_i \in
D$, the process $X$ follows the solution trajectory of the corresponding
initial-value
problem for an exponentially distributed random time, with intensity parameter
$\lambda_i > 0$. Then, a new driving vector
field $u_j$ is selected at random from $D \setminus \{u_i\}$, and $X$ follows
the
corresponding trajectory for another exponentially distributed random time,
with intensity parameter $\lambda_j > 0$. Using exponential times is
required to make the resulting process Markovian. This construction is iterated to
obtain a continuous trajectory on $\R$ that is defined for any time $t \geq 0$
and driven by exactly one of the vector fields in $D$ between any two
switches. We call the intensity parameters $(\lambda_i)_{i \in S}$ switching
rates. For $j \neq i$, let $\lambda_{i,j}$ be the
rate of switching from
$u_i$ to $u_j$. Then, 
\begin{equation*}
 \lambda_i = \sum_{j \neq i} \lambda_{i,j}. 
\end{equation*}

As in~\cite{Bakhtin}, we assume that the exponential times between any two
switches are independent conditioned on the sequence of driving vector
fields, that the switching rate $\lambda_i$ depends only on
the current state $i$ (in particular not on the value of $X$ at
the given time), and that for any two states $i, j \in S$ there is a positive
probability of switching from $i$ to $j$. In many papers on dynamical
systems with random switching, the
switching rates are allowed to depend on the location of the process $X$, and
it is only required that the transition mechanism on $S$ be irreducible (see for
instance~\cite{Faggionato},~\cite{Benaim} and~\cite{Cloez}). We hope to
simplify our exposition by not studying more general classes of
switching systems. 

While $X$ alone is not Markov, we can build a Markov process by
adjoining a second stochastic process $A = (A_t)_{t \geq 0}$ that captures the
current driving vector field at any given time. The
product space $\R \times S$ is then the state space of the 2-component Markov
process $(X,A)$ with Markov semigroup $(\Pp^t)_{t \geq 0}$. We denote the
associated distributions on paths starting at points $(\xi,i)$ by
$\Pp_{\xi,i}$, and the corresponding transition probability measures by
$\Pp^t_{\xi,i}, t \geq 0$. We call $X$ the
continuous and $A$ the discrete component of $(X,A)$. 

Recall that a probability measure $\mu$ on $\R \times S$ is called an
invariant measure of $(\Pp^t)_{t \geq 0}$ if 
\begin{equation}\label{eq:invM}
\mu(E \times \{i\}) = \mu \Pp^t(E \times \{i\}):=\sum_{j \in S} \int_{\R} \Pp^t_{\xi,j}(E \times \{i\}) \
\mu(d\xi \times \{j\}) 
\end{equation}
holds for any Borel set $E \subset \R$, for any $i \in S$ and for any $t \geq
0$. We denote the projection $\mu(\cdot \times \{i\})$ by
$\mu_i$. 

In general, an invariant measure doesn't have to be absolutely
continuous with respect to the product of Lebesgue measure on $\R$ and
counting measure on $S$ (in the sequel, we will simply say ``with
respect to Lebesgue measure''). When it is, the density of the
invariant measure with respect to Lebesgue measure (which is
guaranteed by the Radon--Nikodym theorem) will be called an invariant
density. An invariant density $\rho$ is defined on $\R \times S$, and
we will usually consider the projections $(\rho_i)_{i \in S}$ that are
defined on $\R$ by $\rho_i(\xi) = \rho(\xi,i)$.  In an
abuse of terminology, we refer to these projections as invariant
densities of the invariant measure. These invariant densities are, of
course, elements of $L^1(\R)$ and whenever we state a regularity
property of $\rho_i$, we mean to say that the equivalence class
$\rho_i$ has a representative with this regularity property.

We call a point $\xi \in \R$ {\bf noncritical} if $u_i(\xi) \neq 0$ for all $i
\in S$. We call $\xi \in \R$ {\bf uniformly critical} if $u_i(\xi) = 0$ for all
$i \in S$.
Throughout this paper, we assume that the set of critical points of the vector
fields in $D$ has no accumulation points. If $\xi$ is a
critical point of a vector field $u_i$ for some $i \in S$,
we write that $u_i$ is positive to the right of $\xi$ if there exists an open
interval with left endpoint $\xi$ on which $u_i$ is positive. In this
definition,
``right'' can be replaced with ``left'' and ``positive'' with ``negative''. 

For $i \in S$, we denote the flow function of the vector field $u_i$ by
$\Phi_i$. Due to forward completeness of $u_i$, the flow function is uniquely
defined for any $t \geq 0$ and for any $\eta \in \R$ by 
\begin{align*}
 \frac{d}{dt} \Phi_i^t(\eta) &= u_i(\Phi_i^t(\eta)),  \\
\Phi_i^0(\eta) &= \eta.
\end{align*}
We write $\R_{+}$ to denote the positive real line $(0, \infty)$. For any
vector $\tbf = (t_1, \ldots, t_m) \in \R_{+}^m$ of times between subsequent
switches, and for any corresponding vector $\ibf = (i_1, \ldots, i_m) \in S^m$
of driving states, we define
\begin{equation*}
 \Phi_{\ibf}^{\tbf}(\eta) := \Phi_{i_m}^{t_m}(\Phi_{i_{m-1}}^{t_{m-1}}(\ldots
\Phi_{i_1}^{t_1}(\eta)) \ldots)
\end{equation*}
as the cumulative flow along the trajectories of $u_{i_1}, \ldots, u_{i_m}$
with starting point $\eta \in \R$. 

\section{Main results}  \label{sec:results}

Let $\mu$ be an invariant measure of $(\Pp^t)_{t \geq 0}$ that is absolutely
continuous with
respect to Lebesgue measure. Let $(\rho_i)_{i \in S}$ denote the invariant
densities associated to $\mu$. 
\subsection{Smoothness at noncritical points}

If $n$ is a positive integer, we call a function
$\Cs^n$ on a set $I$ or at a point $\xi$ if the function is $n$ times
continuously differentiable
on $I$ or at $\xi$. Being $\Cs^0$ means being continuous. 

\begin{theorem} \label{thm:smoothness}
Let $\xi \in \R$ be noncritical, and assume that there exist an integer
$n \geq 1$ and a closed interval $I$ containing $\xi$ in its interior on which
all vector fields
in $D$ are $\Cs^{n+1}$. Then, the invariant
densities $(\rho_i)_{i \in S}$ are $\Cs^n$ at $\xi$.  
\end{theorem}

\begin{remark}  \rm  \label{rm:smoothness_2}
The following statement is an immediate consequence of
Theorem~\ref{thm:smoothness}:
If $\xi \in \R$ is noncritical and if all vector fields in $D$ are
$\Cs^{\infty}$ on a closed interval $I$ containing $\xi$ in its interior, then
the invariant
densities are $\Cs^{\infty}$ at $\xi$.
\end{remark}

\begin{remark}  \rm  \label{rm:existence_ac}
According to~\cite[Theorem 1]{Bakhtin}, the following is a sufficient 
condition for absolute continuity (and also for uniqueness) of the invariant
measure on $\R \times S$, provided that an invariant measure exists:

There is a point $\xi \in \R$ that is not uniformly critical and that is
$D$-accessible from any starting point $\eta \in \R$. We say that a point
$\xi$ is $D$-accessible from $\eta$ if for any open
interval $I$ containing $\xi$, there exist a finite sequence of vector fields in
$D$ and a corresponding sequence of switching times such that some point in $I$
can be reached from $\eta$ by following the flows composed of these vector
fields and switching times. 

If the process $X$ is eventually confined to a bounded subset of $\R$,
existence of an invariant measure is guaranteed by the Krylov--Bogoliubov
method. This is for instance the case if each vector field in $D$ has finitely
many critical points, with the smallest critical
point attracting from the left and the largest critical point attracting
from the right.

A sufficient condition for existence and exponential ergodicity that can also be
applied in a noncompact
setting has been provided in~\cite[Assumption 1.8]{Malrieu}:

Let $(\Phi_i)_{i \in S}$ be the flow functions associated to the vector fields
$(u_i)_{i \in S}$, and assume that $\Phi_i^t$ is globally Lipschitz continuous
with
Lipschitz constant $L_i^t$ for any $i \in S$ and for any $t \geq 0$.
Furthermore, suppose that 
\begin{equation*}   \label{f:Wasserstein_curvature}
 \alpha_i := \inf_{t > 0} \Bigl(-\frac{\ln(L_i^t)}{t} \Bigr)
\end{equation*}
is a well-defined real number for any $i \in S$. If one assumes, as we do, that
the switching rates are independent of the position of $X$, the stochastic
process $A$ on $S$ is Markov and has an invariant measure $\nu$. The
condition
\begin{equation}   \label{f:mean_contraction}
 \sum_{i \in S} \nu(\{i\})  \alpha_i > 0
\end{equation}
then implies existence of an invariant measure $\mu$ for $(X,A)$ such that the
distribution of $(X,A)_t$ converges to $\mu$ in a certain Wasserstein distance
as $t$ goes to infinity.

Condition~\ref{f:mean_contraction} can be interpreted as
$(X,A)$ contracting in mean, see~\cite[page 5]{Cloez}.
In~\cite[Assumption 1.3]{Cloez}, the authors state a generalization of
Condition~\ref{f:mean_contraction} for switching between members of a finite
family of Markov processes.
See also~\cite[Sections 3.3 and 3.4]{Yin} for criteria for recurrence and
positive recurrence of the process $(X,A)$.   
\end{remark}

\subsection{Asymptotics at critical points}   \label{subsec:asymptotics}

Let $\xi$ be a critical point of $u_i$ for some $i \in S$, and assume that none
of the other vector fields in $D$
have $\xi$ as a critical point. This assumption is made to simplify the
asymptotic analysis (see~\cite[Section 2]{Balazs}). Without loss of
generality, let $u = u_1$ and let $\xi = 0$. Recall our standing assumption
that for all vector fields $u_j \in
D$, the set of critical points of $u_j$ has no accumulation
point (see Section~\ref{sec:definitions}). Then, there is a $\delta > 0$
such that
none of the vector fields in $D \setminus \{u_1\}$ have a critical point in $[0,
\delta]$ and $u_1$ has no critical point in $(0,\delta]$. To simplify the
analysis further, we assume that there is a constant $a \neq 0$ so that
\begin{equation*}
 u_1(\eta) = -a \eta + O(\eta^2)
\end{equation*}
as $\eta$ approaches $0$ from the right, i.e. $u_1$ behaves almost linearly
near $0$. The constant $a$ can be thought of as the
contraction or expansion coefficient of $u_1$ near $0$.

If $u_1$ was of order $O(\eta^{\alpha})$ for $\alpha < 1$, the vector field
would not be Lipschitz continuous and the resulting ODE could fail to
have unique solutions. If $u_1$ was of order $O(\eta^{\alpha})$ for
$\alpha > 1$, identifying the asymptotically dominant term would become more
complicated. 

Under these assumptions, we study the asymptotic behavior of $\rho_1$ as
$\eta$ approaches $0$ from the right. Due to the symmetric nature of the
problem, there is no need to investigate the case of $\eta$ approaching $0$
from the left separately. 

In Section~\ref{sec:support}, we show that the
support of the measures $(\mu_i)_{i \in S}$
can be represented as a finite union of closed intervals of positive length
(see Corollary~\ref{co:support}). 
Let $\Is$ denote the collection of these intervals. If $\mu$ is ergodic, $\Is$
contains only one interval.  

Exactly one of the following statements holds:

\begin{enumerate}[(A)]
 \item $0$ is the left endpoint of an open interval that does not contain any
points from the support of $(\mu_i)_{i \in S}$. 
\item $0$ is contained in the interior of an interval $I \in \Is$. 
\item $0$ is the left endpoint of an interval $I \in \Is$. 
\end{enumerate}

Although these statements are not formulated in terms of the given vector
fields, it is easy to see which of them holds by using the algorithm at the
beginning of Section~\ref{sec:support}. 
In case A, $\rho_1$ is constantly equal to zero on an open interval with left
endpoint $0$. Cases B and C are more intricate and are dealt with in
Theorems~\ref{thm:critical_points_analytic} and~\ref{thm:critical_points}.
In case C, either $0$ is the right endpoint of an open interval that does
not contain any points from the support of $(\mu_i)_{i \in S}$, or $0$ is the right
endpoint of an interval $J \in \Is$. But if $0$ is both left endpoint of an
interval $I \in \Is$ and right endpoint of an interval $J \in \Is$, it is
uniformly critical (see Section~\ref{sec:support}). Since we assume that $0$
is only critical for $u_1$, this
second scenario cannot occur. 

To illustrate cases B and C, we give two simple examples of PDMPs: one that
corresponds to case B and one that corresponds to case C.  

\begin{example}  \rm \label{ex:case_2}
Let $D$ be the collection of vector fields $u_1(\eta) := -\eta$, $u_2(\eta) :=
1$ and $u_3(\eta):= -1$. At any given time, the process $X$ is either attracted
to the critical point $0$ or moves to the left or to the right at constant
speed. The associated flow functions are $\Phi_1^t(\eta) = \eta e^{-t}$,
$\Phi_2^t(\eta) = \eta + t$ and $\Phi_3^t(\eta) = \eta - t$, with global
Lipschitz constants $L_1^t = e^{-t}$ and $L_2^t = L_3^t = 1$ for any $t \geq
0$. If we define $\alpha_1, \alpha_2, \alpha_3$ as in
Remark~\ref{rm:existence_ac}, we have $\alpha_1 = 1$ and $\alpha_2 =
\alpha_3 = 0$. Since we allow switching from any vector field to any other
vector field, criterion~\eqref{f:mean_contraction} implies existence of an
invariant measure. Theorem 1 in~\cite{Bakhtin} implies that the invariant
measure is unique and absolutely continuous. The
projections $(\mu_i)_{i \in S}$ of this invariant measure
are
supported on $\R$ (see Section~\ref{sec:support}). This is then an example
of case B. 
\end{example}

\begin{example}   \rm  \label{ex:case_3}
Let $D$ be the collection of vector fields $u_1(\eta) := -\eta$ and
$u_2(\eta):= 1 - \eta$. The process $X$ is alternately attracted by $0$ and
$1$, and is eventually confined to the interval $(0,1)$. By the
Krylov--Bogoliubov theorem, $(\Pp^t)_{t \geq 0}$ has an invariant measure. The
invariant measure is unique and absolutely continuous.
The support of
the measures $(\mu_i)_{i \in S}$ is the closed interval
$[0,1]$, so this is an example of case C. 
\end{example}

\subsubsection{Asymptotics for analytic vector fields}

In this subsection, we suppose in addition to the assumptions of
Subsection~\ref{subsec:asymptotics} that all vector fields in $D$ are analytic
in an
open interval around $0$. To state our result on the asymptotically dominant
term of $\rho_1$, we
introduce the function 
\begin{equation*}
 \bar \rho(\eta) := \sum_{i > 1} \lambda_{i,1}  \rho_i(\eta). 
\end{equation*}

\begin{theorem}    \label{thm:critical_points_analytic}
Under the assumptions above, the following statements hold.
\begin{enumerate}
 \item Let $\lambda_1 < a$.
In cases B and C, there is a constant $c > 0$ such that
\begin{equation*}
\rho_1(\eta) = c \eta^{\frac{\lambda_1}{a} -1} +
o(\eta^{\frac{\lambda_1}{a}-1})
\end{equation*}
as $\eta$ approaches $0$ from the right. 
\item Let $\lambda_1 > a > 0$. In case B, 
\begin{equation*}
 \lim_{\eta \downarrow 0} \rho_1(\eta) = \frac{\bar \rho(0)}{\lambda_1 - a} > 0.
\end{equation*}
In case C, there is a constant $c > 0$ such that
\begin{equation*}
 \rho_1(\eta) = c \eta^{\frac{\lambda_1}{a}-1} +
o(\eta^{\frac{\lambda_1}{a}-1})
\end{equation*}
as $\eta$ approaches $0$ from the right. 
\item Let $\lambda_1 = a$. In case B, there is a constant $c > 0$ such
that
\begin{equation*}
 \rho_1(\eta) = -c \ln(\eta) + o(\ln(\eta))
\end{equation*}
as $\eta$ approaches $0$ from the right.
In case C, $\rho_1(\eta)$ converges to a positive constant as $\eta$ approaches
$0$ from the right. 
\end{enumerate}
\end{theorem}
Theorem~\ref{thm:critical_points_analytic} is proved in
Section~\ref{sec:critical_points}. Note that in the critical case $\lambda_1 =
a$, the density $\rho_1$ is
unbounded to the right of $\eta = 0$ in case B and bounded in case C. 

\subsubsection{Asymptotics in the nonanalytic case}

In the absence of analyticity, we can still recover some of the results from
Theorem~\ref{thm:critical_points_analytic}.

\begin{theorem}     \label{thm:critical_points}
Under the assumptions above, without analyticity of the vector fields, the
following statements hold. 
\begin{enumerate}
 \item Let $\lambda_1 < a$.
In cases B and C, there is a constant $c > 0$ such that
\begin{equation*}
\rho_1(\eta) = c \eta^{\frac{\lambda_1}{a} -1} +
o(\eta^{\frac{\lambda_1}{a}-1})
\end{equation*}
as $\eta$ approaches $0$ from the right. 
\item Let $\lambda_1 > a > 0$. In case B, 
\begin{equation*}
 \lim_{\eta \downarrow 0} \rho_1(\eta) = \frac{\bar \rho(0)}{\lambda_1 - a} > 0.
\end{equation*}
In case C,
\begin{equation*}
 \lim_{\eta \downarrow 0} \rho_1(\eta) = \frac{\bar \rho(0)}{\lambda_1 - a} = 0.
\end{equation*}
\item Let $\lambda_1 = a$. In case B, there exist constants $c', c > 0$ such
that
\begin{equation*}
 -c'  \ln(\eta) \leq \rho_1(\eta) \leq -c  \ln(\eta)
\end{equation*}
for $\eta$ sufficiently small. In case C, there is a constant $c > 0$ such
that
\begin{equation*}
 \rho_1(\eta) \leq  -c  \ln(\eta). 
\end{equation*}
\item Let $a < 0$, i.e. $0$ is a repelling critical point of $u_1$. In case B, 
\begin{equation*}
\lim_{\eta \downarrow 0}  \rho_1(\eta) = \frac{\bar \rho(0)}{\lambda_1 - a} > 0.
\end{equation*}
Case C is not possible. 
\end{enumerate}
\end{theorem}
Theorem~\ref{thm:critical_points} is proved in
Section~\ref{sec:critical_points}. Theorem~\ref{thm:critical_points} implies the
following conditions for boundedness of $\rho_1$ to the right of $0$. 

\begin{corollary}     \label{co:boundedness}
\begin{enumerate}
 \item If $\lambda_1 < a$, $\rho_1$ is unbounded to the right of $0$ in cases
B and C.
\item If $\lambda_1 > a > 0$, $\rho_1$ is bounded to the right of $0$ in cases
B and C.
\item If $\lambda_1 = a$, $\rho_1$ is unbounded to the right of $0$ in case
B. In case C, our analysis is inconclusive. 
\end{enumerate}
\end{corollary}

\begin{remark}  \rm \label{rm:boundedness}
The conditions in
Corollary~\ref{co:boundedness} align with
intuition. If $\lambda_1 < a$, the rate of switching away from $u_1$ is lower
than the rate at which $u_1$ contracts to its critical point $0$. In
this case, the rate at which mass accumulates in the vicinity of $0$ is high,
which results in a singularity of the invariant density at $0$. If $\lambda_1 >
 a > 0$, the rate of switching away from $u_1$ is higher than the rate of
contracting to $0$. The rate at which mass accumulates at $0$ is low and
$\rho_1$ is bounded near $0$ (see~\cite[Theorem 1, part c]{Balazs}).
\end{remark}

\section{Integral equations for the densities}  \label{sec:representation}

In this section, we present two integral equations satisfied almost everywhere
by invariant
densities of $(\Pp^t)_{t \geq 0}$. Loosely stated, these equations
illustrate how mass with respect to an invariant density $\rho_i$ accumulates at
a point $\eta \in \R$. At some point in time, there is a switch from a vector
field in $D \setminus \{ u_i\}$ to $u_i$, and the flow associated to $u_i$
transports mass to $\eta$. In a sense, we condition on the time and nature of
this last switch to $u_i$. The family of equations in
Lemma~\ref{lm:integral_equation_1} describe the mass transport for a finite
history of the process. In this case, there is a positive probability of
having no switch. The equation in
Lemma~\ref{lm:integral_equation_2} describes the transport mechanism for an
infinite history. This guarantees that with probability 1, there is at least one
switch.   

Let $\mu$ be an absolutely continuous invariant measure of $(\Pp^t)_{t \geq
0}$, with invariant densities $(\rho_i)_{i \in S}$.
Since we do not assume backward completeness of the vector fields in $D$, we
have to be careful when studying the history of a switching trajectory. It
could happen that the backward flow associated to a vector field goes off to
$-\infty$ or $\infty$ in finite time. For any $\eta \in \R$ and for any $i \in
S$, let $\tau_i(\eta)$ denote the supremum over the set of times $t \geq 0$
for which $t \mapsto \Phi_i^{-t}(\eta)$ is well-defined. With this
definition, 
we introduce the shorthand 
\begin{equation}    \label{f:pushforward}
\Phi_i^t \# h(\eta) := \begin{cases}
	\frac{h(\Phi_i^{-t}(\eta))}{D
\Phi_i^t(\Phi_i^{-t}(\eta))}, & t < \tau_i(\eta)  \\
	0, & t \geq \tau_i(\eta) 
        \end{cases} 
\end{equation}
for the pushforward of the function $h$ under the flow map $\Phi_i^t$. We think
of $h$ as a density function on the real line.   
Note that $D \Phi_i^t >0$ in dimension one, so there is no need for
absolute value in the denominator. Since $u_i$ is assumed to be $\Cs^1$, so is
$\eta \mapsto \Phi_i^t(\eta)$, and the differential $D \Phi_i^t$ is
well-defined.  

Let $L^1_{+}(\R)$ be the set of $L^1$-functions on the real line that have a
nonnegative representative. In other words, $L^1_{+}(\R)$ is the space of
densities for finite measures on $\R$. 

For any $h \in L^1_{+}(\R)$ and for any
$T>0$, define the Perron--Frobenius operators
\begin{equation}  \label{f:Perron_Frobenius_1}
 \bar \Pp^T_i h(\eta):=\frac{1}{T}  \int_0^T e^{-\lambda_i t} 
\,\Phi_i^t \# h(\eta) 
\,dt
\end{equation}
and 
\begin{equation}   \label{f:Perron_Frobenius_2}
\hat \Pp^T_i h(\eta) := \frac{1}{T}
\int_0^T e^{-\lambda_i t}  (T-t) 
\,\Phi_i^t \# h(\eta) \
dt.
\end{equation}

We can now state the truncated version of the integral
equation. 
 
\begin{lemma}  \label{lm:integral_equation_1}
For any $i \in S$ and for any $T > 0$,
\begin{equation*}    
\rho_i \equiv \bar \Pp^T_i \rho_i + \sum_{j \neq i} \lambda_{j,i} 
\hat
\Pp^T_i \rho_j.
\end{equation*} 
\end{lemma}

To state the integral equation over an infinite time horizon, we define the
operators
\begin{equation*}
 \bar \Pp_i h(\eta):= \int_0^{\infty}
e^{-\lambda_i t} 
\,\Phi_i^t \# h(\eta)  \,dt, \quad i \in S
\end{equation*}
for densities $h \in L^1_{+}(\R)$.

\begin{lemma} \label{lm:integral_equation_2}
For any $i \in S$, 
\begin{equation*}
\rho_i = \sum_{j \neq i} \lambda_{j,i}  \bar \Pp_i \rho_j.
\end{equation*} 
\end{lemma}

Lemmas~\ref{lm:integral_equation_1} and~\ref{lm:integral_equation_2} are proved
in Section~\ref{sec:integral_equations}. As will become apparent from these
proofs, the lemmas continue to hold if the state space $\R$ of the continuous
component $X$ is replaced with a finite-dimensional smooth manifold.

\section{Smoothness of
the densities at noncritical points}
\label{sec:smoothness}

In this section, we prove Theorem~\ref{thm:smoothness}, which was stated at the
beginning of Section~\ref{sec:results}. 
By assumption, there exist an integer $n \geq 1$ and a closed interval $I$,
with $\xi$ in its interior,  on which
all vector fields
in $D$ are $\Cs^{n+1}$. Since $\xi$ is noncritical and since for each vector
field in $D$, the set of critical points has no accumulation point, we
may assume without loss of generality that $I$
does not contain any critical points. Let $I_0 \subset I$ be another compact
interval containing $\xi$ in its interior, whose endpoints are a positive
distance
away from the endpoints of $I$. As the trajectories of the $X$-component of the
switching process have bounded speed on compact subsets of $\R$, there is a
small time $T_0 >0$ such that $(\Phi_{\ibf}^{\sbf})^{-1}(\eta) \in I$
for any finite index sequence $\ibf$, any corresponding sequence of
nonnegative switching times $\sbf$ with $l^1$-norm less than or equal to $T_0$
and for any $\eta \in I_0$. 

We define the
integration kernels 
\begin{equation}    \label{f:integration_kernel_1}
\K_i(\zeta,\eta):=\frac{\exp \Bigl(\lambda_i 
\int_{\eta}^{\zeta} \frac{dx}{u_i(x)}
\Bigr)}{u_i(\eta)}
\end{equation}
and
\begin{equation}   \label{f:integration_kernel_2}
\hat \K_i^{T_0}(\zeta,\eta):= \Bigl( T_0 + \int_{\eta}^{\zeta}
\frac{dx}{u_i(x)} \Bigr)  \K_i(\zeta,\eta) 
\end{equation}
for $i \in S$ and $(\zeta,\eta) \in I \times I_0$. 
With these definitions, we have the following representations of $\bar
\Pp^{T_0}_i
\rho_i$ and $\hat \Pp^{T_0}_i \rho_j$. See~\eqref{f:Perron_Frobenius_1}
and~\eqref{f:Perron_Frobenius_2} for the definitions of $\bar \Pp^{T_0}_i$ and
$\hat \Pp^{T_0}_i$.  

\begin{lemma}     \label{lm:integration_kernel_representation}
For any $i \in S$ and for any $\eta \in I_0$, 
\begin{equation}    \label{f:integration_kernel_representation_1}
 \bar \Pp^{T_0}_i \rho_i(\eta) = \frac{1}{T_0} 
\int_{\Phi_i^{-T_0}(\eta)}^{\eta} \rho_i(\zeta)  \K_i(\zeta,\eta) \,d\zeta.
\end{equation}
For any $i, j \in S$, $i \neq j$, and for any $\eta \in I_0$, 
\begin{equation}    \label{f:integration_kernel_representation_2}
 \hat \Pp^{T_0}_i \rho_j(\eta) = \frac{1}{T_0} 
\int_{\Phi_i^{-T_0}(\eta)}^{\eta} \rho_j(\zeta)  \hat
\K_i^{T_0}(\zeta,\eta) \,d\zeta. 
\end{equation}
\end{lemma}
Our definition of $T_0$ implies that the interval
$[\Phi_i^{-T_0}(\eta), \eta]$ (or $[\eta, \Phi_i^{-T_0}(\eta)]$ if $u_i(\xi) <
0$) is contained in $I$, so the integrals on the right are well-defined.
Notice in particular that this reasoning still holds if $u_i$ is not backward
complete. 

\medskip
 
\bpf Fix an $\eta \in I_0$ and recall the definition of $\Phi_i^t \# \rho_i$
in~\eqref{f:pushforward}. Linearity of the Jacobi flow gives 
\begin{equation*}
 D \Phi_i^t(\Phi_i^{-t}(\eta)) =
\frac{u_i(\eta)}{u_i(\Phi_i^{-t}(\eta))},
\end{equation*}
hence
\begin{equation*}
 \Phi_i^t \# \rho_i(\eta) = \rho_i(\Phi_i^{-t}(\eta)) 
\frac{u_i(\Phi_i^{-t}(\eta))}{u_i(\eta)}
\end{equation*}
for any $t \in [0,T_0]$. The change of variables $\zeta =
\Phi_i^{-t}(\eta)$ then
yields~\eqref{f:integration_kernel_representation_1}.
Formula~\eqref{f:integration_kernel_representation_2} is proved similarly.  
\epf 

\bigskip

In~\eqref{f:integration_kernel_representation_1}
and~\eqref{f:integration_kernel_representation_2}, the expressions on the right
still make sense if $\K_i$ and $\hat \K_i^{T_0}$
are replaced with arbitrary kernels on $I \times I_0$. For any such kernel $\K$
and for any $i, j \in S$, set
\begin{equation}    \label{f:integration_kernel_arbitrary}
 \K^{T_0}_i \rho_j(\eta) := \frac{1}{T_0} 
\int_{\Phi_i^{-T_0}(\eta)}^{\eta} \rho_j(\zeta)  \K(\zeta,\eta) \,d\zeta.
\end{equation}

The following lemma addresses regularity of the integration
kernels $(\K_i)_{i \in S}$ and $(\hat \K_i^{T_0})_{i \in S}$.

\begin{lemma} \label{lm:smooth_kernel}
The kernels $(\K_i)_{i \in S}$ and $(\hat
\K_i^{T_0})_{i \in S}$ are $\Cs^{n+1}$ on $I \times I_0$.
\end{lemma}

\bpf This follows from the fact that $u_i$ is $\Cs^{n+1}$ and nonzero on
$I$.  
\epf

\bigskip

The following lemmas illustrate the smoothening effect of the operators $(\bar
\Pp^{T_0}_i)_{i \in S}$ and $(\hat \Pp^{T_0}_i)_{i \in S}$.  We begin by showing
that, away from
critical points, the densities $(\rho_i)_{i \in S}$ are bounded.

\begin{lemma} \label{lm:bounded}
The densities $(\rho_i)_{i \in S}$ are bounded on the interval $I_0$.
\end{lemma}

\bpf Fix an $i \in S$. By Lemma~\ref{lm:integral_equation_1}, it is enough
to show that $\bar \Pp^{T_0}_i \rho_i$ and $(\hat \Pp^{T_0}_i \rho_j)_{j \neq
i}$ are bounded on $I_0$. Since $\K_i$ and $\hat
\K_i^{T_0}$ are $\Cs^1$ on
the compact
set $I \times I_0$ (Lemma~\ref{lm:smooth_kernel}), they are also bounded on $I
\times I_0$ by
constants $k_i$ and $k_{i, T_0}$. For $j \in S$, let $\| \rho_j \|_1$
denote the $L^1$-norm of $\rho_j$ on $\R$. 
Using~\eqref{f:integration_kernel_representation_1},
\begin{equation*}
 \bar \Pp^{T_0}_i \rho_i(\eta) \leq \frac{k_i}{T_0} 
\int_{\Phi_i^{-T_0}(\eta)}^{\eta} \rho_i(\zeta) \,d\zeta \leq \frac{k_i
 \| \rho_i \|_1}{T_0}
\end{equation*}
for any $\eta \in I_0$.

And using~\eqref{f:integration_kernel_representation_2},
\begin{equation*}
\hat \Pp^{T_0}_i \rho_j(\eta) \leq \frac{k_{i, T_0}}{T_0}
 \int_{\Phi_i^{-T_0}(\eta)}^{\eta} \rho_j(\zeta) \,d\zeta \leq
\frac{k_{i, T_0}  \| \rho_j \|_1}{T_0}
\end{equation*}
for any $j \neq i$, $\eta \in I_0$.
\epf

\bigskip

\begin{remark}   \rm
In the proof of Lemma~\ref{lm:bounded}, we did not use any concrete
information about $\K_i$ or $\hat \K_i^{T_0}$ other than boundedness on $I
\times
I_0$. The result
still holds if $\K_i$ and $\hat \K_i^{T_0}$ are replaced with arbitrary kernels
that
are bounded on $I
\times I_0$. Furthermore, the time $T_0$ can be replaced with any time $T \in
(0,T_0)$.  
\end{remark}

The following corollary will be useful in Section~\ref{sec:critical_points}
when we derive asymptotics for invariant densities at critical points. 

\begin{corollary}  \label{co:bounded}
Let $i \in S$ and assume that $\xi \in \R$ is not a critical point of $u_i$.
Then, there is a compact interval $I$ with $\xi$ in its
interior, such that $\rho_i$ is bounded on $I$. 
\end{corollary}

In Lemma~\ref{lm:bounded}, we assumed that $\xi$ is noncritical. Here, the
point $\xi$ may be critical for some of the vector fields in $D$, just not for
the particular vector field $u_i$ whose corresponding density function we are
interested in. 

\medskip

\bpf Since $u_i(\xi) \neq 0$ and since the set of critical points of $u_i$ has
no accumulation points, there is a compact interval $I$ that has $\xi$ in
its interior, but does
not contain any critical points of $u_i$. Let $I_0 \subset I$ be another
compact interval with $\xi$ in its interior such that the endpoints of $I_0$
are a positive distance away from the endpoints of $I$. Choose $T>0$ so small
that $\Phi_i^{-t}(\eta) \in I$ for any $\eta \in I_0$ and for any $t \in [0,T]$.
Define the kernels $\K_i$ and
$\hat \K_i^T$ according to~\eqref{f:integration_kernel_1}
and~\eqref{f:integration_kernel_2}. These kernels are bounded on $I \times
I_0$, and
we can repeat the proof of Lemma~\ref{lm:bounded} to finish the
argument. 
\epf

\bigskip

Let $I_1 \subset I_0$ be
a compact interval that contains $\xi$ in its interior
and whose endpoints are a positive distance away from the endpoints of $I_0$.
Let $T_1 \in (0,T_0]$ be so small that $(\Phi_{\ibf}^{\sbf})^{-1}(\eta) \in
I_0$ for any index sequence $\ibf$, any corresponding sequence of
nonnegative
switching times $\sbf$ with
$l^1$-norm less than or equal to $T_1$, and for any $\eta \in I_1$.

\begin{lemma} \label{lm:Lipschitz}
The densities $(\rho_i)_{i \in S}$ are
Lipschitz continuous
on $I_1$.  
\end{lemma}

\bpf Fix an $i \in S$. By Lemma~\ref{lm:integral_equation_1}, it is enough
to show that $\bar \Pp^{T_1}_i \rho_i$ and $(\hat \Pp^{T_1}_i \rho_j)_{j \neq
i}$ are Lipschitz continuous on $I_1$. By Lemma~\ref{lm:bounded}, $\rho_i$ is
bounded on $I_0$ by some constant $r_i$. Let $L$ be a Lipschitz constant of
$\K_i$ on $I \times I_0$ and let $\tilde L$ be a Lipschitz constant of the flow
function $\Phi_i$ on $[-T_1,0] \times I_1$. The constant $k_i$ is defined as in
the proof of Lemma~\ref{lm:bounded} and $k_{i, T_1}$ is defined in analogy
to $k_{i, T_0}$. Fix two points $\eta,
\vartheta \in I_1$. As $\Phi_i^{-T_1}(\eta)$ and $\Phi_t^{-T_1}(\vartheta)$ are
both contained in $I_0$, we obtain the estimate
\begin{align*}
& \lvert \bar \Pp^{T_1}_i \rho_i(\eta) - \bar \Pp^{T_1}_i \rho_i(\vartheta)
\rvert \\
=& \frac{1}{T_1}  \Bigl \lvert \int_{\Phi_i^{-T_1}(\eta)}^{\eta}
\rho_i(\zeta)  \K_i(\zeta,\eta) \,d\zeta -
\int_{\Phi_i^{-T_1}(\vartheta)}^{\vartheta}
\rho_i(\zeta)  \K_i(\zeta,\vartheta) \,d\zeta \Bigr \rvert \\
\leq& \frac{1}{T_1}  \Bigl( \Bigl \lvert
\int_{\Phi_i^{-T_1}(\eta)}^{\Phi_i^{-T_1}(\vartheta)} \rho_i(\zeta) 
\K_i(\zeta,\eta) \,d\zeta \Bigr \rvert + \Bigl \lvert \int_{\eta}^{\vartheta}
\rho_i(\zeta)  \K_i(\zeta,\vartheta) \,d\zeta \Bigr \rvert  \\
& + \Bigl \lvert \int_{\Phi_i^{-T_1}(\vartheta)}^{\eta} \rho_i(\zeta) 
(\K_i(\zeta,\eta) - \K_i(\zeta,\vartheta)) \,d\zeta \Bigr \rvert \Bigr) \\
\leq& \| \vartheta - \eta \|  \frac{1}{T_1}  (r_i k_i  (1
+ \tilde L) + L  \| \rho_i \|_1).
\end{align*}

Let $\hat L$ be a Lipschitz constant of $\hat \K_i^{T_1}$ on $I \times I_0$. For
a fixed $j \neq i$, the density $\rho_j$ is bounded on $I_0$ by a constant
$r_j$, and 
\begin{equation}    \label{f:Lipschitz}
\lvert \hat \Pp^{T_1}_i \rho_j(\eta) - \hat \Pp^{T_1}_i \rho_j(\vartheta)
\rvert \leq \lvert \vartheta -\eta
\rvert  \frac{1}{T_1}  (r_j k_{i, T_1} 
(1 +
\tilde L) + \hat L  \| \rho_j \|_1).
\end{equation}
\epf

\bigskip

\begin{remark} \rm \label{rm:Lipschitz}
Lemma~\ref{lm:Lipschitz} continues to hold if $\K_i$ and $\hat \K_i^{T_1}$ are
replaced
with arbitrary
kernels that are Lipschitz continuous on $I \times I_0$ and if $T_1$ is
replaced with an arbitrary time
$T \in (0,T_1)$. 
\end{remark}

\begin{remark}  \rm \label{rm:continuity}
Lemma~\ref{lm:Lipschitz} implies the following: If an open interval $I$ does
not contain
 any critical points, then all invariant densities $\rho_i$ are Lipschitz
continuous
 on $I$. Slightly modifying the proof of Lemma~\ref{lm:Lipschitz}, one can show
a related statement:
 
 If an open interval $I$ does not contain any critical points of a particular
vector field $u_i$
 (but possibly critical points of other vector fields), the invariant density
$\rho_i$ is continuous
 on $I$. 
 
Notice that we can only guarantee the continuity of
$\rho_i$, not its Lipschitz continuity. Since we allow for critical 
 points of the other vector fields $(u_j)_{j \neq i}$ on $I$, we also can no longer
ascertain boundedness of the 
 corresponding densities $(\rho_j)_{j \neq i}$. Instead of~\eqref{f:Lipschitz},
we obtain the weaker estimate 
\begin{align*}
\lvert \hat \Pp^{T_1}_i \rho_j(\eta) - \hat \Pp^{T_1}_i
\rho_j(\vartheta) \rvert \leq& \frac{k_{i, T_1}}{T_1}
 \Bigl( \Bigl \lvert 
 \int_{\Phi_i^{-T_1}(\vartheta)}^{\Phi_i^{-T_1}(\eta)} \rho_j(\zeta) \,d\zeta
\Bigr \rvert 
 + \Bigl \lvert \int_{\vartheta}^{\eta} \rho_j(\zeta) \,d\zeta \Bigr \rvert
\Bigr) \\
 & + \frac{\| \rho_j \|_1}{T_1}  \hat L  \lvert \vartheta - \eta
\rvert.
\end{align*}
\end{remark}

Lemma~\ref{lm:differentiability} illustrates the actual smoothening mechanism. 

\begin{lemma} \label{lm:differentiability}
For any integer $k \in \{0, \ldots, n-1\}$, the following statement holds:

If the densities $(\rho_i)_{i \in S}$ are
$\Cs^k$ on a compact interval $\tilde I \subset I_1$ that contains $\xi$ in its
interior, there exist a compact interval $\tilde I' \subset \tilde I$
with $\xi$ in its interior and a time $T \in (0, T_1]$ such that for any
$\Cs^{k+2}$-kernel $\K$ on $I \times I_0$, the
functions $(\K^T_i \rho_j)_{i,j \in S}$ are $\Cs^{k+1}$ on $\tilde I'$.
\end{lemma}
Recall that we defined $\K^T_i \rho_j$
in~\eqref{f:integration_kernel_arbitrary}. 

\medskip

\bpf We prove Lemma~\ref{lm:differentiability} by induction on $k$. In the base
case, assume that the densities $(\rho_i)_{i \in S}$ are continuous on
$\tilde I \subset I_1$.  Let $\tilde I' \subset
\tilde I$ be a compact interval that contains $\xi$ in its interior and whose
endpoints are a positive distance away from the endpoints of $\tilde I$. Let $T
\in (0, T_1]$ be so small that $(\Phi_{\ibf}^{\sbf})^{-1}(\eta) \in \tilde
I$ for any index sequence $\ibf$, any corresponding sequences of nonnegative
switching times $\sbf$ with $l^1$-norm less than or equal to $T$, and for any
$\eta \in \tilde I'$. For any $\Cs^2$-kernel $\K$ on $I \times I_0$, for any
$\eta \in
\tilde I'$ and for any $i, j \in S$,
\begin{align}  
\frac{d}{d\eta} \K^{T}_i \rho_j(\eta) =& \frac{1}{T}  \Bigl(\rho_j(\eta)
 \K(\eta,\eta) -
\rho_j(\Phi_i^{-T}(\eta))  \K(\Phi_i^{-T}(\eta),\eta) 
\frac{d}{d\eta}\Phi_i^{-T}(\eta) \Bigr)  \notag \\
& + \frac{1}{T}  \int_{\Phi_i^{-T}(\eta)}^{\eta} \rho_j(\zeta) 
\partial_2
\K(\zeta,\eta) \,d\zeta   \notag \\
=& \frac{1}{T}  \Bigl(\rho_j(\eta)  \K(\eta,\eta) -
\rho_j(\Phi_i^{-T}(\eta))  \K(\Phi_i^{-T}(\eta),\eta) 
\frac{d}{d\eta}\, \Phi_i^{-T}(\eta) \Bigr)  \label{f:differentiability} \\  
& + (\partial_2 \K)^{T}_i
\rho_j(\eta).  \notag
\end{align}
Here, $\partial_2 \K$ denotes the partial derivative of $\K$ with respect
to its second component. Since $\rho_j$ is assumed to be $\Cs^0$ on
$\tilde I$, since $\K$ is $\Cs^2$ on $I \times I_0$ and since $u_i$ is
$\Cs^{n+1}$ on
$I$, the first term in~\eqref{f:differentiability} is $\Cs^0$ on $\tilde I'$. 

It remains to show that $(\partial_2 \K)^{T}_i \rho_j$ is $\Cs^0$, but this
follows along the lines of Lemma~\ref{lm:Lipschitz},
keeping in mind that
$\partial_2
\K$ is $\Cs^1$ and hence Lipschitz continuous on $I \times I_0$
and that $T \leq T_1$ (see
Remark~\ref{rm:Lipschitz}). 
Since $\tfrac{d}{d\eta} \K^{T}_i
\rho_j(\eta)$ is $\Cs^0$ on $\tilde I'$, it follows that $\K^{T}_i
\rho_j$
is $\Cs^1$ on $\tilde I'$. This completes the base case. 


In the induction step, let $k$ be a fixed integer in $\{1, \ldots, n-1\}$ and
assume that the statement holds for $k-1$.
Assume that the densities $(\rho_i)_{i \in S}$ are $\Cs^k$ on $\tilde I \subset
I_1$. The densities
$(\rho_i)_{i \in S}$ are then also $\Cs^{k-1}$ on $\tilde I$. By induction
hypothesis, there exist a compact interval $\tilde I' \subset
\tilde I$ with $\xi$ in its interior and a time $T \in (0,T_1]$ such that for
any $\Cs^{k+1}$-kernel $\K$ on $I \times I_0$, the functions $(\K_i^T
\rho_j)_{i, j \in S}$ are $\Cs^k$ on $\tilde I'$. Without loss of generality, we
can assume that the endpoints of $\tilde I'$ are a positive distance away from
the endpoints of $\tilde I$ and that $T$ is so small
that $(\Phi_{\ibf}^{\sbf})^{-1}(\eta) \in \tilde I$ for any index sequence
$\ibf$, any corresponding sequence of nonnegative switching times $\sbf$ with
$l^1$-norm less than or equal to $T$, and for any $\eta \in \tilde I'$. Let $\K$
be a $\Cs^{k+2}$-kernel on $I \times I_0$. Then,~\eqref{f:differentiability}
holds for any $\eta \in \tilde I'$
and for any $i, j \in S$. 

Since $\rho_j$ is by assumption $\Cs^k$ on
$\tilde I$, since $\K$ is $\Cs^{k+2}$ on $I \times I_0$ and since $u_i$ is
$\Cs^{n+1}$ on
$I$, the first term in~\eqref{f:differentiability} is $\Cs^k$ on $\tilde I'$.
In addition, $\partial_2 \K$ is a
$\Cs^{k+1}$-kernel on $I \times I_0$. By induction hypothesis, $(\partial_2
\K)^T_i \rho_j$
is $\Cs^k$ on $\tilde I'$, so $\tfrac{d}{d\eta} \K^{T}_i
\rho_j(\eta)$ is $\Cs^k$ on $\tilde I'$. From this, it follows that
$\K^{T}_i \rho_j$
is $\Cs^{k+1}$ on $\tilde I'$. 
\epf

\bigskip

\bpf [Proof of Theorem~\ref{thm:smoothness}]
In order to prove Theorem~\ref{thm:smoothness}, it suffices to show the
following statement: For any integer $k \in \{0, \ldots, n\}$, there is a compact interval
$I_{k+1}$ with $\xi$ in its interior such that the densities $(\rho_i)_{i \in
S}$ are $\Cs^k$ on $I_{k+1}$. 

We prove this statement by induction on $k$. By Lemma~\ref{lm:Lipschitz},
the densities $(\rho_i)_{i \in S}$ are Lipschitz continuous on $I_1$. This
takes care of the base case.  
In the induction step, let $k$ be an integer in $\{1, \ldots, n\}$ and assume
that the densities $(\rho_i)_{i \in S}$ are $\Cs^{k-1}$ on a compact interval
$I_k \subset I_1$ with $\xi$ in its interior.
By Lemma~\ref{lm:differentiability}, there exist a compact interval $I_{k+1}
\subset I_k$ with $\xi$ in its interior and a time $T \in (0, T_1]$ such that
for any $\Cs^{k+1}$-kernel $\K$ on $I \times I_0$, the functions $(\K_i^T
\rho_j)_{i, j \in S}$ are $\Cs^k$ on $I_{k+1}$. 

Fixing an $i \in S$, Lemma~\ref{lm:differentiability} applied to the integration
kernel $\K_i$ yields that $\bar \Pp^{T}_i \rho_i$ is $\Cs^k$ on
$I_{k+1}$. Applying Lemma~\ref{lm:differentiability} to $\hat \K_i^T$
yields that for any $j \neq i$, $\hat \Pp^{T}_i \rho_j$ is $\Cs^k$ on $I_{k+1}$.
By Lemma~\ref{lm:integral_equation_1}, $\rho_i$ is $\Cs^k$ on $I_{k+1}$.
\epf


\section{The support of invariant measures on the real line}
\label{sec:support}

Let $\mu$ be an invariant measure of the Markov semigroup $(\Pp^t)_{t \geq
0}$. In this section, we describe the support of the measures $(\mu_i)_{i
\in S}$, which are measures on the real line. 
Recall that a point $\xi \in \R$ lies in the support of $\mu_i$ if and
only if $\mu_i(U) >0$ for any open neighborhood $U$
of $\xi$. 

We say that a point $\xi \in \R$ is $D$-reachable from a point $\eta \in \R$ if
there exist a finite index sequence $\ibf$ and a corresponding sequence
of
nonnegative
switching times $\tbf$ such that
\begin{equation*}
 \Phi_{\ibf}^{\tbf}(\eta) = \xi. 
\end{equation*}
For any $\xi \in \R$, we define $L(\xi)$ as the set of
points that are $D$-reachable from $\xi$. 

We call a point $\xi \in \R$ $D$-accessible from $\eta \in \R$ if for
any
open neighborhood $U$ of $\xi$ there exist a finite index sequence $\ibf$ and a
corresponding sequence of nonnegative switching times $\tbf$ such that
\begin{equation*}
 \Phi_{\ibf}^{\tbf}(\eta) \in U. 
\end{equation*}
Let $L$ denote the set of points on the real line that are $D$-accessible from
any point in $\R$. 

\subsection{Minimal invariant sets}

A nonempty set $I \subset \R$ is called invariant if 
\begin{equation*}
 \Phi_{\ibf}^{\tbf}(\xi) \in I
\end{equation*}
for any $\xi \in I$, any finite index sequence $\ibf$ and any
corresponding sequence of nonnegative switching times $\tbf$. A minimal
invariant set is an invariant
set for which any nonempty strict subset is not invariant. Alternatively, a
minimal invariant set is a nonempty set $I$ with the property that
\begin{equation}   \label{f:minimal_invariant_set}
 L(\xi) = I
\end{equation}
for any $\xi \in I$. 

The following algorithm yields exactly the minimal invariant sets of our
switching system. 
\begin{enumerate}
 \item Mark $-\infty$ with the label ``l'' and mark $+\infty$ with the
label
``r''. 
\item Mark those critical points where all vector fields in $D$ are
nonnegative with an ``l'' and mark those critical points where all vector
fields in $D$ are nonpositive with an ``r''. If a critical point has both
labels ``l'' and ``r'', it is uniformly critical. All uniformly critical
points form minimal invariant sets.  
\item Consider all points, including $-\infty$, with the label ``l''. This
includes those points that carry both labels. As $+\infty$ doesn't have label
``l'', each of these points has a closest labeled point to its right. If this
point has label ``r'', the open, possibly infinite, interval with the
``l''-labeled and the ``r''-labeled points as its endpoints is a candidate for
a minimal invariant set. It is indeed a minimal invariant set if and only if it
contains two not necessarily distinct points $\xi$ and $\eta$ for which there
are vector fields $u, v \in
D$ with $u(\xi)>0$ and $v(\eta) <0$. 
\end{enumerate}

\begin{proposition}  \label{prop:minimal_invariant_sets}
The algorithm above characterizes the minimal invariant sets of the
switching system completely. Minimal invariant sets are thus either open
intervals or point sets with exactly one element.  
\end{proposition}

\bpf We first show that any set identified by the algorithm as a minimal
invariant set is indeed a minimal invariant set. Let $\Sc$ be a set identified
by the algorithm as a minimal invariant set. Then,
either $\Sc = \{\xi\}$, where $\xi$ is a uniformly critical point, or $\Sc$
is an open interval $(l,r)$, where $l < r$ are elements of the extended real
line such that 
\begin{enumerate}[(a)]
 \item $l= -\infty$ or $u_i(l) \geq 0$ for any $i \in S$
 \item $r= \infty$ or $u_i(r) \leq 0$ for any $i \in S$
 \item for any critical point $\xi$ in $(l,r)$ there exist indices $i,j \in S$
with $u_i(\xi) < 0 < u_j(\xi)$
 \item if there are no critical points in $(l,r)$, there are at least points
$\xi, \eta \in (l,r)$ and indices $i, j \in S$ with
$u_i(\xi) < 0 < u_j(\eta)$.
\end{enumerate}
If $\Sc = \{\xi\}$, it is clear that $\Sc$ is a minimal invariant set: The only
strict subset of $\Sc$ is the empty set, and $\Sc$ is invariant because $\xi$
is uniformly critical. 

If $\Sc = (l,r)$, no switching trajectory starting in $\Sc$ can get to the left
of $l$ or to the right of $r$. This is obvious if $l = -\infty$ or $r =
\infty$. If $l$ or $r$ are finite, it is guaranteed by Conditions a and b,
respectively.  
Hence, $\Sc$ is invariant. Next, we show that $\Sc$ is also minimal. Assume
that $\Sc$ is not minimal. Then, there is a nonempty strict subset $\Rc$
of $\Sc$ that is invariant. In addition, there is a point $\zeta \in \Sc$
with $u_i(\zeta) \leq 0$ for any $i \in S$. To see this, fix a point
$\eta \in \Sc \setminus \Rc$ and a point
$\xi \in \Rc$. We can assume without loss of generality that $\eta > \xi$.
Since $\Rc$ is invariant, $\eta$ is not $D$-reachable from $\xi$. Thus, there
is a point $\zeta \in [\xi,\eta]$ with $u_i(\zeta) \leq 0$ for any $i \in S$.

In light of Condition c, $\zeta$ is not critical. On the other hand, Condition
d implies that there exist a $\tilde
\zeta \in \Sc$ and a $j \in S$
with $u_j(\tilde \zeta) > 0$. 
Assume without loss of generality that $\tilde \zeta > \zeta$, and let 
\begin{equation*}
 \hat \zeta := \sup \{ \theta \in [\zeta,\tilde \zeta]: \,u_i(\theta) < 0 \
\forall i \in S\}. 
\end{equation*}
The point $\hat \zeta$ is a critical point in $\Sc$ with $u_i(\hat
\zeta) \leq 0$ for any $i \in S$. This violates Condition c. 

\medskip

Conversely, let $I$ be a minimal invariant set. We need to show that the
algorithm correctly identifies $I$ as a minimal invariant set. Due to the
minimality
assumption, $I$ is an interval. If $I$ contains exactly one point, this
point is uniformly critical, for otherwise, $I$ would not be invariant. 

If $I$ has at least two elements, it is an interval with distinct
endpoints. We now
show that if an
endpoint of $I$ is finite, it must be a critical point.
Let $\xi$ be a finite endpoint of $I$, say its left endpoint, and assume that
$\xi$ is noncritical. Since $I$ is invariant, $u_i(\xi) > 0$ for
any $i \in S$. By continuity of the vector fields, there is an $\epsilon > 0$
such that $u_i(\eta) > 0$ for any $i \in S$ and for any $\eta \in [\xi, \xi +
\epsilon]$. By choosing $\epsilon$ sufficiently small, we can then ensure
that $I \setminus [\xi,\xi + \epsilon]$ is a nonempty strict subset of $I$ that
is invariant. This contradicts the minimality assumption on $I$. 

Invariance of $I$ also implies that the endpoints of $I$ are not
$D$-reachable from a starting
point in the interior of $I$. Hence, $I$ is an
open interval $(l,r)$, where $l$ and $r$ may be finite
or infinite. 

It remains to show that Conditions c and d are satisfied. Let $\xi \in I$ be a
critical point. If $u_i(\xi) \geq 0$ for any $i \in S$, the interval
$(\xi,r) \subset (l,r)$ is invariant as well which is a contradiction. Similarly,
$u_i(\xi)$ cannot be nonpositive for all $i \in S$, so we can find $i, j \in
S$ with $u_i(\xi) < 0 < u_j(\xi)$. To show that Condition d holds, assume that
$u_i(\eta) \geq 0$ for all $\eta \in I$ and for all $i \in S$. For $\xi \in I$,
the interval $(\xi,r)$ is invariant, which contradicts the
minimality assumption. 
\epf 

\bigskip

\begin{proposition}   \label{prop:disjoint}
Minimal invariant sets are pairwise disjoint.
\end{proposition}

\bpf
Let $I$ and $J$ be minimal invariant sets with $I \cap J \neq \emptyset$. As
the intersection of invariant sets, $I \cap J$ is invariant. Since $I$ and $J$
are minimal, it follows that $I = I \cap J = J$. 
\epf

\bigskip

\bigskip

\subsection{How minimal invariant sets relate to the support of invariant
measures}

We begin by relating invariant measures of the global dynamics on $\R
\times S$ to invariant measures of the switching dynamics confined to a minimal
invariant set. 

Let $I \subset \R$ be a minimal invariant set. On $I \times S$, we define the
semigroup $(\pp^t)_{t \geq 0}$ by 
\begin{equation*}
 \pp^t_{\xi,i}(E \times \{j\}) := \Pp^t_{\xi,i}(E \times \{j\})
\end{equation*}
for any $(\xi,i) \in I \times S$ and for any set $E$ in the Borel
$\sigma$-algebra on $I$. Hence, $(\pp^t)_{t \geq 0}$ can be thought of as the
restriction of $(\Pp^t)_{t \geq 0}$ to $I \times S$. It is well-defined because
$I$ is invariant. 

\begin{proposition}   \label{prop:restricted_semigroup}
There is a one-to-one correspondence between invariant measures of $(\pp^t)_{t
\geq 0}$ and those invariant measures of $(\Pp^t)_{t \geq 0}$ that assign mass
$1$ to $I \times S$. 
\end{proposition}

This is easy to see. We omit the proof. 

\medskip

Next, we show that the support of the
measure $\mu(\cdot \times S)$ does
not contain
points outside of the closure of minimal invariant sets.

\begin{proposition}  \label{prop:support_in_mis}
Let $\xi \in \R$ be a point that is not contained in the closure of a minimal
invariant set. Then,
$\xi$ is not contained in the support of $\mu(\cdot
\times S)$.  
\end{proposition}

\bpf [Sketch of proof] We record two statements without proofs.

First, there exist an open neighborhood $U$ of
$\xi$, an open set
$V \subset \R$, a time $T>0$, a positive integer $m$, an index sequence $\ibf =
(i_1, \ldots, i_m)$
of length $m$ and an open subset $\Delta$ of the simplex 
\begin{equation*}
 \Bigl\{\sbf \in (0, \infty)^{m-1}: \,\sum_{l=1}^{m-1} s_l < T \Bigr\}
\end{equation*}
 such that $\Phi_{\ibf}^{(\sbf,T-\sum_{l=1}^{m-1} s_l)}(\eta) \in V$ for any
$\eta \in U, \sbf \in \Delta$, and such that $U \cap L(\vartheta) =
\emptyset$ for any $\vartheta \in V$.

Second, there are constants $\epsilon', c' > 0$ and an
open set $U' \subset U$ such that
\begin{equation*}
 \inf_{\eta \in U'; i,j \in S} \Pp_{\eta,i}^{\epsilon'}(U \times \{j\})
\geq c'.
\end{equation*}
The second statement is shown in~\cite{Benaim}. It is an immediate consequence
of the fact that the speed of the process $X$ is bounded on bounded sets. 
 
To derive a contradiction, we assume that $\xi$ belongs to the support of
$\mu(\cdot \times S)$. As $U'$
is an open neighborhood of $\xi$,  
\begin{equation*}
c:= \mu(U' \times S) > 0.
\end{equation*}
Therefore, 
\begin{align*}
\sum_{i \in S} \int_{U'} \Pp^{T+\epsilon'}_{\theta,i}(V \times S) \
\mu_i(d\theta) =& \sum_{i \in S} \int_{U'} \sum_{l \in S} \int_{\R}
\Pp^T_{\eta,l}(V \times S)
 \Pp^{\epsilon'}_{\theta,i}(d\eta \times \{l\}) \,\mu_i(d\theta) \\
\geq& \sum_{i \in S} \int_{U'} \int_U \Pp^T_{\eta,i_1}(V \times S) 
\Pp^{\epsilon'}_{\theta,i}(d\eta \times \{i_1\}) \,\mu_i(d\theta)
\\
\geq& c  c'  \inf_{\eta \in U} \Pp^T_{\eta,i_1}(V \times S).
\end{align*}
Next, we show that $\inf_{\eta \in U} \Pp^T_{\eta,i_1}(V \times S) > 0$. Fix a
point $\eta \in U$ and let $C_{\ibf}$ denote the event that the driving
vector fields up to time $T$ appear in the order given by $\ibf$. Let
$\Pp_{\ibf}$ be the probability that the first $m$ driving vector fields appear
in the order given by $\ibf$, conditioned on $u_{i_1}$ being the first driving
vector field. Similarly to Lemma~\ref{lm:identity} in
Section~\ref{sec:integral_equations}, we have  
\begin{align}
\Pp^T_{\eta,i_1}(V \times S) \geq& \Pp_{\eta,i_1}(X_T \in V, C_{\ibf}) \notag \\
\geq& \Pp_{\ibf}  \int_{\Delta}
\prod_{l=1}^{m-1} \lambda_{i_l} 
e^{-\lambda_{i_l} s_l}  e^{-\lambda_{i_m} (T - (s_1 + \ldots +
s_{m-1}))}
d\sbf.  \label{f:support_in_mis}
\end{align}
The term in~\eqref{f:support_in_mis} is strictly positive and does not depend
on $\eta$. 
We conclude that  
\begin{equation*}
a:= \sum_{i \in S} \int_{U'} \Pp^{T+\epsilon'}_{\theta,i}(V \times S) \
\mu_i(d\theta) > 0.
\end{equation*}

Hence, there is a positive integer $N$ with $N  a >1$. Let
$\mu \Pp$ denote the distribution of the Markov process $(X,A)$ with initial
distribution $\mu$. For $0 \leq k \leq N-1$, define the event 
\begin{equation*}
E_k:=\{X_{k  (T+\epsilon')} \in U', X_{(k+1) 
(T+\epsilon')} \in V, X_{j  (T+\epsilon')} \in U'^{c} \ \text{for} \ k+2
\leq j \leq N\}.
\end{equation*}
Since the sets $(E_k)_{0 \leq k \leq N-1}$ are pairwise disjoint,
\begin{equation*}
\mu \Pp(X_{N  (T +\epsilon')} \in U'^{c}) \geq \sum_{k=0}^{N-1}
\mu \Pp(E_k).
\end{equation*}
Since $U'$ cannot be reached from any point in $V$, we have 
\begin{equation*}
 E_k = \{X_{k  (T+\epsilon')} \in U', X_{(k+1)  (T+\epsilon')} \in
V\}.
\end{equation*}
Then, for $0 \leq k \leq N-1$,
\begin{equation*}
\mu \Pp(E_k) = \sum_{i \in S} \int_{\R} \sum_{j \in S} \int_{U'}
\Pp^{T+\epsilon'}_{\eta,j}(V \times S)  \Pp^{k 
(T+\epsilon')}_{\theta,i}(d\eta \times \{j\}) \,\mu_i(d\theta) = a 
\end{equation*}
because $\mu$ is an invariant measure. We infer that
\begin{equation*}
\mu \Pp(X_{N  (T+\epsilon')} \in U'^{c}) \geq N  a > 1,
\end{equation*}
which is impossible. Hence, $\xi$ is not contained in the support of $\mu(\cdot
\times S)$.  
\epf

\bigskip

In Proposition~\ref{prop:restricted_semigroup}, we saw that invariant
measures on minimal invariant sets correspond to invariant measures on $\R$
that are supported on a minimal invariant set. In the following proposition,
we show uniqueness of the invariant measure on a given minimal invariant set. 

\begin{proposition}     \label{prop:unique_im}
Any minimal invariant set admits at most one invariant measure.
\end{proposition}

\bpf Let $I$ be a minimal invariant set. If $I = \{\xi\}$ for some uniformly
critical point $\xi$, uniqueness of the invariant measure is clear. 

If $I$ is an open interval, it does not contain any uniformly critical points
by Proposition~\ref{prop:minimal_invariant_sets}. By the alternative
characterization of minimal invariant sets in~\eqref{f:minimal_invariant_set},
$I = L(\eta)$ for any $\eta \in I$. Thus, any point in $I$ is $D$-reachable
from all starting points in $I$. By~\cite[Theorem 1]{Bakhtin}, this
implies uniqueness of the invariant measure of the restricted semigroup
$(\pp^t)_{t \geq 0}$. 
\epf

\bigskip

Now, assume that the invariant measure $\mu$ is ergodic. If $I$ is a
minimal invariant
set, ergodicity of $\mu$ implies that $\mu(I \times S)$ is either $0$ or $1$.
It is then natural to ask whether we can assign a unique minimal invariant set
$I$ to $\mu$ for which $\mu(I \times S) = 1$. The following proposition shows
that this can be done. 

\begin{proposition}   \label{prop:unique_mis}
If $\mu$ is ergodic, there is a unique minimal invariant
set $I$ with $\mu(I \times S) = 1$.
\end{proposition}

\bpf Let us begin by showing that such a minimal invariant set exists.
Since $\mu$
is ergodic, it is enough to show that $\mu(I \times S) > 0$ for some minimal
invariant set $I$. 

We denote the set of points not contained in the closure of a
minimal invariant set by $\Tc$. According to
Proposition~\ref{prop:support_in_mis}, the intersection of $\Tc$ and of the
support of $\mu(\cdot \times S)$ is empty, so there is a point $\xi \in
\Tc^c$ that also lies in the support of $\mu(\cdot \times S)$. As $\xi \in
\Tc^c$, there is a minimal invariant set $I$ whose closure contains $\xi$.
We distinguish between several cases.

First, assume that $I = \{\xi\}$. Then, $\xi$ is uniformly critical and may or
may not be an endpoint of one or two additional minimal invariant sets. If
there are no minimal invariant sets adjacent to $\{\xi\}$, we can find an open
neighborhood $U$ of $\xi$ such that $U \setminus \{\xi\} \subset \Tc$. Since the
complement of the support of $\mu(\cdot \times S)$ has measure $0$, it follows
that $\mu(U \setminus \{\xi\} \times S) = 0$. Therefore, $\mu(\{\xi\} \times S)
> 0$. 

If there is at least one open minimal invariant set adjacent to $\{\xi\}$, we
have $\mu(\{\xi\} \times S) > 0$, or at least one of the adjacent minimal
invariant sets has strictly positive $\mu(\cdot \times S)$-measure. 

Now, assume that $I = (l,r)$. If $\xi \in I$, it is immediate from the
definition of the support that $\mu(I \times S) > 0$. If $\xi$ is an
endpoint of $I$, assume without loss of generality that $\xi = l$.
We have
already dealt with the case where $\xi$ is uniformly critical. If $\xi$ is
critical but not uniformly critical, we still have $\mu(\{\xi\} \times S) > 0$
or $\mu(I \times S) >0$ or $\mu(J \times S) > 0$, provided that $J$ is an open
minimal invariant set with $\xi$ as its right endpoint. We only need to
exclude the case that $\mu(\{\xi\} \times S) > 0$. This can be done similarly to
the proof of Proposition~\ref{prop:support_in_mis}.

\medskip

Uniqueness of the minimal invariant set follows from
Proposition~\ref{prop:disjoint}. 
\epf

\bigskip

\begin{proposition}    \label{prop:support}
If $\mu$ is ergodic, there is a unique minimal
invariant
set $I$ such that the support of the measures $(\mu_i)_{i
\in
S}$
equals the closure of $I$. 
\end{proposition}

\bpf 
Let $I$ be the unique minimal
invariant set for which $\mu(I \times S) = 1$ and whose existence is postulated
in Proposition~\ref{prop:unique_mis}. By
characterization~\eqref{f:minimal_invariant_set} of minimal invariant sets,
every point in $I$ is $D$-reachable from any other point in $I$. By~\cite[Lemma
6]{Bakhtin}, $I$ is then contained in the support of $(\mu_i)_{i \in S}$. With
$\mu(I \times S) = 1$, the statement follows.   
\epf

\bigskip

\begin{corollary}   \label{co:support}
Let $\mu$ be an invariant measure, not necessarily ergodic. Then, there exist
minimal invariant sets $I_1, \ldots, I_N$ such that the support of $\mu_j$
equals the closure of $\bigcup_{i=1}^N I_i$ for any $j \in S$.
If $\mu$ is absolutely continuous, each of the minimal invariant sets
$I_i$ is an open interval. 
\end{corollary}

\bpf
This follows from Proposition~\ref{prop:support}, the Birkhoff Ergodic
Theorem and the 
Ergodic Decomposition Theorem. See~\cite{Hairer:ergodicity-lectures}
for discussion adapted precisely to this setting
and \cite{sinai} for more general considerations.
\epf

\bigskip

\section{Asymptotics for the densities at critical points} 
\label{sec:critical_points}

\subsection{Asymptotic analysis for nonanalytic vector fields} 
\label{subsec:critical_points}

In this subsection, we prove Theorem~\ref{thm:critical_points}.
In both cases B and C, there is an open
interval $I$ with left endpoint $0$ such that
\begin{equation*}
 \rho_i(\eta) > 0
\end{equation*}
for any $\eta \in I$ and for any $i \in S$. This is a consequence of the
following lemma.

\begin{lemma}  \label{lm:positive_densities}
Let $I$ be an open interval that is contained in the support of
$(\mu_i)_{i \in S}$. If the vector field $u_i$ does not have any critical
points in $I$, then $\rho_i(\eta) > 0$ for any $\eta \in I$. 
\end{lemma}

\bpf Fix a point $\eta \in I$. Let $\tilde I$ be a closed subinterval of $I$,
with
$\eta$
contained
in the interior of $\tilde I$. Let $T > 0$ be so small
that $\Phi_{\ibf}^{\sbf}(\eta) \in \tilde I$ for any finite index
sequence $\ibf$ and any corresponding sequence of switching times $\sbf$
with $\| \sbf \|_1 \leq T$. 

Since $u_i$ does not have any critical points in $I$, the function
\begin{equation*}
 \zeta \mapsto \exp \Bigl( -\lambda_i  \int_{\zeta}^{\eta}
\frac{dx}{u_i(x)} \Bigr)
\end{equation*}
is bounded below on $[\Phi_i^{-T}(\eta),\eta]$ by a constant $c > 0$. Using
Lemma~\ref{lm:integral_equation_1}
and~\eqref{f:integration_kernel_representation_1}, we obtain the estimate
\begin{align}
 \rho_i(\eta) \geq& \frac{1}{u_i(\eta)} 
\frac{1}{T} 
\int_{\Phi_i^{-T}(\eta)}^{\eta} \rho_i(\zeta)  \exp \Bigl(
-\lambda_i  \int_{\zeta}^{\eta} \frac{dx}{u_i(x)} \Bigr) \
d\zeta  \notag \\
\geq& \frac{c}{T}  \frac{1}{\lvert u_i(\eta) \rvert} 
\Bigl
\lvert \int_{\Phi_i^{-T}(\eta)}^{\eta} \rho_i(\zeta) \,d\zeta \Bigr
\rvert \notag \\
=& \frac{c}{T}  \frac{1}{\lvert u_i(\eta) \rvert} 
\mu_i((\Phi_i^{-T}(\eta),\eta)) > 0.  \label{f:positive_densities}
\end{align}
For~\eqref{f:positive_densities}, we used that $(\Phi_i^{-T}(\eta), \eta)$ is
contained in the support of $\mu_i$. 
\epf

\bigskip

Let $\delta > 0$ be so
small that none of the vector fields $(u_i)_{i > 1}$ have a critical point in
$[0,\delta]$, and that
$u_1$ has no critical point in $(0, \delta]$. Let $a > 0$. The vector field
$u_1$ is then strictly negative on $(0,\delta]$.  

For $\eta \in (0,\delta)$, define $\vartheta := \lim_{t \to \tau_1(\eta)}
\Phi_1^{-t}(\eta)$, where $\tau_1(\eta)$ was introduced in
Section~\ref{sec:representation}. This limit is independent of the concrete
choice of $\eta$.
By Lemma~\ref{lm:positive_densities}, there is a
constant $c > 0$ such that $\bar \rho(\eta) \geq c$
for any $\eta \in [\tfrac{\delta}{2}, \delta]$. In case B, we can even assume
that $\bar \rho(\eta) \geq c$ for any $\eta \in
[0,\delta]$. And by Remark~\ref{rm:continuity}, $\bar \rho$ is 
continuous on $[0,\delta]$, which implies that $\bar \rho$ is
bounded from above on $[0,\delta]$ by some constant $\bar \rho_{\infty}$.

Set 
\begin{equation*}
 r(\eta) := -\frac{1}{u_1(\eta)} - \frac{1}{a \eta}, \quad \eta \in
(0,\vartheta).
\end{equation*} 
It is not hard to see that $r(\eta)$ is bounded on $(0,\delta]$ by a
constant $r_{\infty} > 0$.
Furthermore, as $u_1 < 0$ on $(0,\vartheta)$, we have $r(\eta) \geq -\frac{1}{a
\eta}$ for any $\eta \in (0,\vartheta)$. 
  
For $\eta, \zeta \in [0,\vartheta]$, define
\begin{equation*}
 E(\eta, \zeta) := \exp \Bigl( -\lambda_1  \int_{\eta}^{\zeta} r(x) \,dx
\Bigr).
\end{equation*}

\begin{lemma}     \label{lm:integrable_function}
 The function $\zeta \mapsto \zeta^{-\frac{\lambda_1}{a}}  \bar
\rho(\zeta)  E(\eta,\zeta)$ is integrable on $(\delta,\vartheta)$ for any
$\eta \in [0,\delta]$. 
\end{lemma}

\bpf For $\eta \in [0,\delta]$ and $\zeta \in (\delta,\vartheta)$,  
\begin{align*}
\zeta^{-\frac{\lambda_1}{a}}  \bar \rho(\zeta)  E(\eta,\zeta) 
=& \zeta^{-\frac{\lambda_1}{a}}  \bar \rho(\zeta)  E(\eta,\delta)
 E(\delta,\zeta) \\
\leq& \zeta^{-\frac{\lambda_1}{a}} 
\bar \rho(\zeta)  e^{\lambda_1 \delta r_{\infty}}  \exp \Bigl(
\lambda_1  \int_{\delta}^{\zeta} \frac{dx}{ax} \Bigr) \\
=& \bar \rho(\zeta)  e^{\lambda_1 \delta r_{\infty}} 
\delta^{-\frac{\lambda_1}{a}}. 
\end{align*}
The fact that $\bar \rho$ is integrable implies the statement. 
\epf

\bigskip

In analogy to Lemma~\ref{lm:integration_kernel_representation}, we have the
following representation for $\rho_1$.

\begin{lemma}   \label{lm:representation_infinite_time}
For any $\eta \in (0,\delta)$,
\begin{equation*}
 \rho_1(\eta) = \Bigl( \frac{\eta^{\frac{\lambda_1}{a} - 1}}{a} + r(\eta)
\eta^{\frac{\lambda_1}{a}} \Bigr)  \int_{\eta}^{\vartheta}
\zeta^{-\frac{\lambda_1}{a}}  \bar \rho(\zeta)  E(\eta,\zeta) \
d\zeta. 
\end{equation*} 
\end{lemma}

\bpf Fix an $\eta \in (0,\delta)$. Using Lemma~\ref{lm:integral_equation_2} and
the change of variables $\zeta = \Phi_1^{-t}(\eta)$, we obtain 
\begin{equation}    \label{f:representation_infinite_time}
 \rho_1(\eta) = -\frac{1}{u_1(\eta)}  \int_{\eta}^{\vartheta} \bar
\rho(\zeta)  \exp \Bigl(\lambda_1  \int_{\eta}^{\zeta}
\frac{dx}{u_1(x)}
\Bigr) \,d\zeta. 
\end{equation}
Since 
\begin{equation*}
 \exp \Bigl( \lambda_1  \int_{\eta}^{\zeta} \frac{dx}{u_1(x)} \Bigr) 
 = \exp \Bigl( -\lambda_1  \int_{\eta}^{\zeta} \frac{dx}{ax} \Bigr)
 E(\eta,\zeta) = \eta^{\frac{\lambda_1}{a}} 
\zeta^{-\frac{\lambda_1}{a}} 
E(\eta,\zeta)
\end{equation*}
for any $\zeta \in (\eta,\vartheta)$, and since $\zeta \mapsto
\zeta^{-\frac{\lambda_1}{a}}  \bar
\rho(\zeta)  E(\eta,\zeta)$ is integrable by
Lemma~\ref{lm:integrable_function}, the statement follows.
\epf

\bigskip

\bpf  [Proof of Theorem~\ref{thm:critical_points}] Fix an $\eta \in (0,\delta)$.
Throughout the proof, we work with the formula for $\rho_1$
provided in Lemma~\ref{lm:representation_infinite_time}.

\medskip

First, let $\lambda_1 < a$. Observe that $\zeta \mapsto
\zeta^{-\frac{\lambda_1}{a}}  \bar \rho(\zeta)  E(0,\zeta)$ is
integrable on $(0,\delta)$ because 
\begin{equation*}
 \zeta^{-\frac{\lambda_1}{a}}  \bar \rho(\zeta) 
E(0,\zeta) \leq \zeta^{-\frac{\lambda_1}{a}}  \bar \rho_{\infty}
e^{\lambda_1 \delta r_{\infty}}. 
\end{equation*}
Together with Lemma~\ref{lm:integrable_function}, we see that this function is
integrable on $(0,\vartheta)$, which implies that
\begin{equation*}
 \lim_{\eta \downarrow 0} \Bigl( \int_{\eta}^{\vartheta}
\zeta^{-\frac{\lambda_1}{a}}  \bar \rho(\zeta)  E(\eta,\zeta) \
d\zeta \Bigr) = \int_0^{\vartheta} \zeta^{-\frac{\lambda_1}{a}}  \bar
\rho(\zeta)  E(0,\zeta) \,d\zeta < \infty
\end{equation*}
by dominated convergence. In addition,
\begin{align*}
\int_0^{\vartheta} \zeta^{-\frac{\lambda_1}{a}}  \bar \rho(\zeta)  
E(0,\zeta) \,d\zeta &\geq \int_{\frac{\delta}{2}}^{\delta}
\zeta^{-\frac{\lambda_1}{a}}  \bar
\rho(\zeta) 
 E(0,\zeta) \,d\zeta \\
&\geq \frac{\delta^{1-\frac{\lambda_1}{a}}}{2}  c  e^{-\lambda_1
\delta r_{\infty}} > 0.
\end{align*}
And since $r(\eta)$ is bounded on $(0,\delta)$, $\lim_{\eta \downarrow 0}
(r(\eta)  \eta^{\frac{\lambda_1}{a}}) = 0$. Part 1 of
Theorem~\ref{thm:critical_points} follows then from
Lemma~\ref{lm:representation_infinite_time}.  

\medskip

Now, let $\lambda_1 > a > 0$. In case B, $\bar \rho(0) > 0$
by Lemma~\ref{lm:positive_densities}. In case C, $\bar \rho(0) = 0$
because $0$ is the right endpoint of an open interval that does not contain any
points from the support of $(\mu_i)_{i \in S}$. 

Since $\lambda_1 > a > 0$, there is a small $\alpha > 0$ such that
\begin{equation*}
\frac{\lambda_1}{a}  (1-\alpha) > 1. 
\end{equation*}
Let $\eta \in (0,\delta)$ so that $\eta < \eta^{\alpha} < \delta$. Then,
\begin{align}
 \rho_1(\eta) =& \Bigl( \frac{\eta^{\frac{\lambda_1}{a} - 1}}{a} + r(\eta)
\eta^{\frac{\lambda_1}{a}} \Bigr)  \int_{\eta}^{\eta^{\alpha}}
\zeta^{-\frac{\lambda_1}{a}}  \bar \rho(\zeta)  E(\eta,\zeta) \
d\zeta \label{f:critical_points_lambda_greater_a_1} \\
& + \Bigl( \frac{\eta^{\frac{\lambda_1}{a} - 1}}{a} + r(\eta)
\eta^{\frac{\lambda_1}{a}} \Bigr)  \int_{\eta^{\alpha}}^{\vartheta}
\zeta^{-\frac{\lambda_1}{a}}  \bar
\rho(\zeta)  E(\eta,\zeta) \,d\zeta.
\label{f:critical_points_lambda_greater_a_2}
\end{align}
The term in~\eqref{f:critical_points_lambda_greater_a_2} is bounded from above
by
\begin{align*}
 \Bigl( \frac{\eta^{\frac{\lambda_1}{a} - 1}}{a} &+ r_{\infty}
\eta^{\frac{\lambda_1}{a}} \Bigr)  \Bigl( \int_{\eta^{\alpha}}^{\delta}
\zeta^{-\frac{\lambda_1}{a}}  \bar
\rho(\zeta)  E(\eta,\zeta) \,d\zeta + \int_{\delta}^{\vartheta}
\zeta^{-\frac{\lambda_1}{a}}  \bar \rho(\zeta)  E(\eta,\zeta) \
d\zeta \Bigr)\\
\leq& \Bigl( \frac{\eta^{\frac{\lambda_1}{a} - 1}}{a} + r_{\infty}
\eta^{\frac{\lambda_1}{a}} \Bigr)  ( \eta^{-\frac{\lambda_1}{a}
\alpha}  e^{\lambda_1 \delta r_{\infty}}
 \| \bar \rho \|_{1} + \delta^{-\frac{\lambda_1}{a}}  e^{\lambda_1
\delta r_{\infty}}  \| \bar \rho \|_1 ) \\
=& e^{\lambda_1 \delta r_{\infty}}  \| \bar \rho \|_1 
 \Bigl( \frac{\eta^{\frac{\lambda_1}{a}
 (1-\alpha) -1}}{a} + r_{\infty} \eta^{\frac{\lambda_1}{a} 
(1-\alpha)} + \frac{\delta^{-\frac{\lambda_1}{a}}}{a}
\eta^{\frac{\lambda_1}{a}-1} + r_{\infty} \delta^{-\frac{\lambda_1}{a}}
\eta^{\frac{\lambda_1}{a}} \Bigr),
\end{align*}
which converges to $0$ as $\eta$ approaches $0$ from the right. 

Since $\eta^{\alpha} < \delta$, the function $\zeta \mapsto \bar \rho(\zeta)
 E(\eta,\zeta)$ is continuous on $[\eta,\eta^{\alpha}]$. By the
mean-value theorem for integration, there exists $\zeta_{\eta} \in
(\eta,\eta^{\alpha})$ such that the term to the right of the equality sign
in~\eqref{f:critical_points_lambda_greater_a_1} equals
\begin{align*}
 \Bigl( \frac{\eta^{\frac{\lambda_1}{a} - 1}}{a} &+ r(\eta)
\eta^{\frac{\lambda_1}{a}} \Bigr)  \int_{\eta}^{\eta^{\alpha}}
\zeta^{-\frac{\lambda_1}{a}} \,d\zeta  \bar
\rho(\zeta_{\eta})  E(\eta,\zeta_{\eta}) \\
=& \Bigl( \frac{1}{a}  (1 - \eta^{(1 - \alpha) 
(\frac{\lambda_1}{a} - 1)}) + r(\eta)  (\eta - \eta^{\alpha +
(1-\alpha)  \frac{\lambda_1}{a}}) \Bigr)  \frac{a}{\lambda_1 - a}
 \bar
\rho(\zeta_{\eta})  E(\eta,\zeta_{\eta}).
\end{align*}
Since $\zeta_{\eta} \in (\eta,\eta^{\alpha})$ for any $\eta$, it is clear that
$\lim_{\eta \downarrow 0} \zeta_{\eta} = 0$. Continuity of $\bar \rho$ at $\eta
= 0$ and integrability of $r(x)$ on $(0,\delta)$ then imply that
\begin{equation*}
\lim_{\eta \downarrow 0} \Bigl( \bar \rho(\zeta_{\eta}) 
E(\eta,\zeta_{\eta}) 
\frac{a}{\lambda_1 - a}
\Bigr) = a  \frac{\bar \rho(0)}{\lambda_1 - a}.
\end{equation*}
Furthermore,
\begin{equation*}
 \lim_{\eta \downarrow 0} \Bigl( \frac{1}{a}  (1 - \eta^{(1 - \alpha)
(\frac{\lambda_1}{a} - 1)}) \Bigr) = \frac{1}{a}. 
\end{equation*}
Finally, for small $\eta > 0$, we have $\eta > \eta^{\alpha + (1-\alpha) 
\frac{\lambda_1}{a}}$.
It follows that
\begin{equation*}
\lvert r(\eta) \rvert  (\eta - \eta^{\alpha + (1-\alpha) 
\frac{\lambda_1}{a}})
\leq r_{\infty}  (\eta - \eta^{\alpha + (1-\alpha)
 \frac{\lambda_1}{a}}),
\end{equation*}
which converges to $0$ as $\eta$ approaches $0$ from the right. This completes
the proof of part 2 of Theorem~\ref{thm:critical_points}. 

\medskip

Next, assume that $\lambda_1 = a$. For $\eta \in
(0,\delta)$, 
\begin{align}
 \rho_1(\eta) =& \Bigl( \frac{1}{a} + r(\eta)
\eta \Bigr)  \int_{\eta}^{\delta}
\zeta^{-1}  \bar \rho(\zeta)  E(\eta,\zeta) \
d\zeta \notag \\
& + \Bigl( \frac{1}{a} + r(\eta)
\eta \Bigr)  \int_{\delta}^{\vartheta}
\zeta^{-1}  \bar
\rho(\zeta)  E(\eta,\zeta) \,d\zeta. 
\label{f:critical_points_lambda_equals_a}
\end{align}
By Lemma~\ref{lm:integrable_function}, the term
in~\eqref{f:critical_points_lambda_equals_a} is bounded on $(0,\delta)$. 
In case B, $c \leq \bar \rho(\eta) \leq \bar \rho_{\infty}$ for any $\eta \in
[0,\delta]$. Therefore,
\begin{align*}
& - c  e^{-\lambda_1 \delta r_{\infty}}  \ln(\eta) + c 
e^{-\lambda_1 \delta r_{\infty}}
 \ln(\delta) \\
\leq& \int_{\eta}^{\delta} \zeta^{-1}  \bar \rho(\zeta) 
E(\eta,\zeta) \,d\zeta \\
\leq& -\bar \rho_{\infty}  e^{\lambda_1 \delta r_{\infty}}  \ln(\eta) 
+ \bar \rho_{\infty}  e^{\lambda_1 \delta r_{\infty}}  \ln(\delta)
\end{align*}
for $\eta \in (0,\delta)$. As  
\begin{equation*}
\lim_{\eta \downarrow 0} \Bigl( \frac{1}{a} + r(\eta) \eta \Bigr) =
\frac{1}{a}, 
\end{equation*}
this establishes part 3 of Theorem~\ref{thm:critical_points} for case B. 
In case C, we only have $\bar \rho(\eta) \leq \bar
\rho_{\infty}$, which is why we obtain a weaker statement. The proof of part 4
is similar to the proof of part 2 and we omit it.  
 
\epf

\bigskip

\subsection{Asymptotic analysis for analytic vector fields} 
\label{ssec:asymptotics_analytic}

In this subsection, we prove Theorem~\ref{thm:critical_points_analytic}. The
ensuing paragraph follows~\cite{Balazs}. 

For any $i \in S$,
we introduce
the probability flux 
\begin{equation*}
 \varphi_i(\eta) = \rho_i(\eta)  u_i(\eta).
\end{equation*}
The vector of probability fluxes $(\varphi_1(\eta), \ldots,
\varphi_n(\eta))^T$ is denoted by $\varphi(\eta)$. As in
Subsection~\ref{subsec:critical_points}, we let $\delta > 0$ be so small that
the vector fields $(u_i)_{i > 1}$ have no critical point in $[0,\delta]$ and
$u_1$ has no critical point in $(0,\delta]$. 

Since the invariant densities $(\rho_i)_{i \in S}$ are $\Cs^1$ on $(0,\delta)$,
they satisfy the Fokker--Planck
equations
\begin{equation}    \label{f:Fokker_Planck}
 \rho_i'(\eta) u_i(\eta) + \rho_i(\eta) u_i'(\eta) = -\lambda_i \rho_i(\eta) +
\sum_{l \neq i} \lambda_{l,i} \rho_l(\eta), \quad i \in S,
\end{equation}
on $(0,\delta)$, see~\cite{Faggionato}. Written in terms of the probability
fluxes, \eqref{f:Fokker_Planck} becomes
\begin{equation}    \label{f:Fokker_Planck_flux}
 \varphi_i'(\eta) = -\frac{\lambda_i}{u_i(\eta)}  \varphi_i(\eta) +
\sum_{l \neq i} \frac{\lambda_{l,i}}{u_l(\eta)}  \varphi_l(\eta), \quad i
\in S.
\end{equation}
In Appendix~\ref{sec:Fokker_Planck}, we show how
Equation~\eqref{f:representation_infinite_time} can be derived directly from the
Fokker--Planck equations if the invariant densities are $\Cs^1$. 

Our approach is to derive the asymptotically dominant term for the probability
flux $\varphi_1$, which will then immediately give the asymptotically dominant
term for $\rho_1$. We begin by showing that $\lim_{\eta \downarrow 0}
\varphi_1(\eta) = 0$.

\begin{lemma}    \label{lm:vanishing_flux}
 We have $\lim_{\eta \downarrow 0} \varphi_1(\eta) = 0$.
\end{lemma}

\bpf By Remark~\ref{rm:continuity}, the limit $\lim_{\eta \downarrow 0}
\varphi_i(\eta)$ exists for any $i > 1$.

It is an
easy corollary of~\eqref{f:Fokker_Planck_flux} that
\begin{equation*}
 \sum_{i \in S} \varphi_i'(\eta) = 0
\end{equation*}
for any $\eta \in (0,\delta)$. Thus, the sum of all probability fluxes is equal
to a constant $k$ on this interval. Since
\begin{equation*}
 \varphi_1(\eta) = k-\sum_{i > 1} \varphi_i(\eta)
\end{equation*}
for any $\eta \in (0,\delta)$, the limit $l := \lim_{\eta \downarrow 0}
\varphi_1(\eta)$ exists as well. It remains to show that $l=0$.

To obtain a contradiction, assume that $l \neq 0$. Then, there is no loss
of generality in assuming that 
\begin{equation*}
 \lvert \varphi_1(\eta) \rvert \geq \frac{\lvert l \rvert}{2}
\end{equation*}
for any $\eta \in (0,\delta)$. Since  $u_1(\eta) = -a \eta + o(\eta)$ as $\eta$
approaches $0$ from the right, we may also assume that
\begin{equation*}
\biggl \lvert \frac{u_1(\eta)}{\eta} \biggr \rvert \leq 2 \lvert a
\rvert, \quad \eta \in (0,\delta). 
\end{equation*}
But this yields
\begin{equation*}
 \int_0^{\delta} \rho_1(\eta) \,d\eta = \int_0^{\delta}
\frac{\lvert \varphi_1(\eta) \rvert}{\lvert u_1(\eta) \rvert} \,d\eta \geq
\frac{\lvert l \rvert}{4 \lvert a \rvert}  \int_0^{\delta} \frac{d
\eta}{\eta} = \infty, 
\end{equation*}
which contradicts the fact that $\rho_1$ is integrable. 
\epf

\bigskip

\begin{corollary}   \label{co:vanishing_flux}
 In case C from Section~\ref{subsec:asymptotics}, $\lim_{\eta \downarrow 0} \varphi(\eta) = 0$.
\end{corollary}

\bpf In case C, the invariant densities $(\rho_i)_{i \in S}$ vanish to the left
of $0$. By Remark~\ref{rm:continuity}, the densities $(\rho_i)_{i > 1}$ are
continuous at $0$, which implies that $\lim_{\eta \downarrow 0} \rho_i(\eta) =
0$ for any $i > 1$. Hence, $\lim_{\eta \downarrow 0} \varphi_i(\eta) = 0$ for
any $i > 1$, and $\lim_{\eta \downarrow 0} \varphi_1(\eta) = 0$ by
Lemma~\ref{lm:vanishing_flux}. 
\epf

\bigskip

We introduce the matrix of switching rates
\begin{equation*}
 \Lambda:= 
 \begin{pmatrix}
  -\lambda_1 & \lambda_{2,1} & \cdots & \lambda_{n,1} \\
  \lambda_{1,2} & -\lambda_2 & \cdots & \lambda_{n,2} \\
  \vdots  & \vdots  & \ddots & \vdots  \\
  \lambda_{1,n} & \lambda_{2,n} & \cdots & -\lambda_n
 \end{pmatrix},
\end{equation*}
and let $U(\eta)$ be the diagonal matrix with diagonal entries
$\tfrac{1}{u_1(\eta)}, \ldots, \tfrac{1}{u_n(\eta)}$. 

For a fixed $\epsilon \in (0,\delta)$, we
consider the initial-value problem
\begin{align}  
 	\phi'(\eta) &= \Lambda  U(\eta)  \phi(\eta), \label{f:ivp_1}
\\
	\phi(\epsilon) &= \varphi(\epsilon),  \notag
\end{align}
whose unique solution is $\varphi(\eta)$. Initial-value problem~\eqref{f:ivp_1}
can be written equivalently as
\begin{align}
      \phi'(\eta) &= \frac{1}{\eta} B(\eta)  \phi(\eta), \label{f:ivp_2}\\
      \phi(\epsilon) &= \varphi(\epsilon). \notag
 \end{align}
Here,
\begin{equation*}    
 B(\eta) := \Lambda  \tilde U(\eta),
\end{equation*}
where $\tilde U(\eta)$ is the diagonal matrix with diagonal entries
$\tfrac{\eta}{u_1(\eta)}, \ldots, \tfrac{\eta}{u_n(\eta)}$. Note that 
$B(\eta)$ is analytic at $\eta = 0$. This follows from the fact that the
diagonal entries of
$\tilde U(\eta)$ are analytic at $\eta = 0$, which is easily derived from
analyticity of the vector fields. The linear system~\eqref{f:ivp_2} then has a
so-called regular singular point at $\eta = 0$ (see~\cite[Section
3.11]{Taylor_ode}). 

Since $B(\eta)$ is analytic at $\eta = 0$, there exist a $\rho \in
(0,\delta)$ and a sequence of matrices $(B_k)_{k \geq 0}$ such that
\begin{equation}  \label{f:expansion}
 B(\eta) = \sum_{k=0}^{\infty} \eta^k  B_k
\end{equation}
for any $\eta \in (-\rho,\rho)$. There is no loss of generality in assuming
that $\rho = \delta$. 

Since $u_1(\eta) = -a \eta + O(\eta^2)$, and since $u_i(\eta) \neq 0$
for any $i > 1$, the matrix $B_0$ in~\eqref{f:expansion} has the
form
\begin{equation*}
B_0 = \begin{pmatrix}
  \frac{\lambda_1}{a} & 0 & \cdots & 0 \\
  -\frac{\lambda_{1,2}}{a} & 0 & \cdots & 0 \\
  \vdots  & \vdots  & \ddots & \vdots  \\
  -\frac{\lambda_{1,n}}{a} & 0 & \cdots & 0
 \end{pmatrix}. 
\end{equation*}
It is easy to give a complete description of the eigenvalues and corresponding
eigenspaces of $B_0$.

\begin{lemma}  \label{lm:eigenvalues_B_0}
The matrix $B_0$ has eigenvalues $\tfrac{\lambda_1}{a}$ and $0$. The
eigenspace \\
corresponding to $\tfrac{\lambda_1}{a}$ is spanned by the vector
$\lambda:=(\lambda_1, -\lambda_{1,2}, -\lambda_{1,3}, \ldots,
-\lambda_{1,n})^T$. The eigenspace corresponding to $0$ is the orthogonal
complement to the span of $\{(1,0,\ldots,0)^T\}$.
\end{lemma}
We omit the proof of Lemma~\ref{lm:eigenvalues_B_0}.

\medskip

At this point, we need to
distinguish between two cases. First, assume that $\tfrac{\lambda_1}{a}$ is not
an integer. Such a condition is
sometimes referred to as
a nonresonance condition. The following statement is then a reformulation
of~\cite[Proposition 11.2]{Taylor_ode}.

\begin{lemma}  \label{lm:matrix_de}
There is a function
\begin{equation}     \label{f:expansion_of_V}
 V(\eta) = \mathbbm{1} + \sum_{k=1}^{\infty} \eta^k  V_k
\end{equation}
that satisfies the normal equation 
\begin{equation}     \label{f:normal_equation}
 \eta V'(\eta) = B(\eta) V(\eta) - V(\eta) B_0, \quad \eta \in (0,\delta)
\end{equation}
and for which 
\begin{equation*}
 \varphi(\eta) = V(\eta)  \exp \Bigl( \ln \Bigl( \frac{\eta}{\epsilon}
\Bigr) B_0
\Bigr) V(\epsilon)^{-1} \varphi(\epsilon), \quad \eta \in (0,\delta).
\end{equation*}
\end{lemma}

\medskip

Now, we consider the resonance case, i.e we assume that $\tfrac{\lambda_1}{a}$
is a positive integer. In this case, we may not be able to construct a solution
of the form~\eqref{f:expansion_of_V} to~\eqref{f:normal_equation}. Instead, we
consider the modified version 
\begin{equation}    \label{f:modified_normal_equation}
 \eta V'(\eta) = B(\eta) V(\eta) - V(\eta) (B_0 + \eta^{\frac{\lambda_1}{a}}
Y), 
\end{equation}
where $Y$ is a matrix satisfying 
\begin{equation}    \label{f:defining_equation_Y}
 B_0 Y = Y \Bigl(B_0 + \frac{\lambda_1}{a} \mathbbm{1} \Bigr).
\end{equation}
In this setting, we have the following reformulation
of~\cite[Proposition 11.5]{Taylor_ode}.

\begin{lemma}   \label{lm:integer_eigenvalue}
There exist a function $V(\eta)$ of the form~\eqref{f:expansion_of_V} and a
matrix $Y$ satisfying~\eqref{f:defining_equation_Y} such that $V(\eta)$
satisfies~\eqref{f:modified_normal_equation} with $Y$ and
\begin{equation*}
 \varphi(\eta) = V(\eta)  \exp \Bigl( \ln \Bigl( \frac{\eta}{\epsilon}
\Bigr) B_0
\Bigr)  \exp \Bigl( \ln \Bigl( \frac{\eta}{\epsilon} \Bigr) Y
\Bigr) V(\epsilon)^{-1} \varphi(\epsilon), \quad \eta \in (0,\delta).
\end{equation*}
\end{lemma}

\medskip

\bpf[Proof of Theorem~\ref{thm:critical_points_analytic}] Comparing
Theorems~\ref{thm:critical_points_analytic}
and~\ref{thm:critical_points}, we see that we only need to show part 2 for
case C and part 3 for both cases. 

Let $\nu \in \R$ and let $\tilde y \in \R^n$ with first component equal to $0$
such that 
\begin{equation*}
 V(\epsilon)^{-1} \varphi(\epsilon) = \nu \lambda + \tilde y, 
\end{equation*}
where $\lambda$ was defined in Lemma~\ref{lm:eigenvalues_B_0}.

In the nonresonance case, Lemma~\ref{lm:eigenvalues_B_0} implies that
\begin{align}
 \exp \Bigl( \ln \Bigl( \frac{\eta}{\epsilon} \Bigr) B_0 \Bigr)
V(\epsilon)^{-1} \varphi(\epsilon) &= \sum_{k=0}^{\infty} \frac{1}{k!} 
\Bigl( \ln \Bigl(
\frac{\eta}{\epsilon} \Bigr) \Bigr)^k (\nu B_0^k \lambda +
B_0^k \tilde y) \notag \\
&= \tilde y + \nu \lambda + \sum_{k=1}^{\infty} \frac{1}{k!}  \Bigl( \ln
\Bigl( \frac{\eta}{\epsilon} \Bigr) \Bigr)^k \nu \Bigl(
\frac{\lambda_1}{a} \Bigr)^k \lambda \notag \\
&= \tilde y + \nu  \exp \Bigl( \frac{\lambda_1}{a}  \ln \Bigl(
\frac{\eta}{\epsilon} \Bigr) \Bigr) \lambda \notag \\
&= \tilde y + \nu \epsilon^{-\frac{\lambda_1}{a}} \eta^{\frac{\lambda_1}{a}}
\lambda,  \label{f:critical_points_analytic_4}
\end{align}
so
\begin{equation}     \label{f:critical_points_analytic_1}
 \varphi(\eta) = \Bigl( \mathbbm{1} + \sum_{k=1}^{\infty} \eta^k V_k \Bigr)
 (
\tilde y + \nu \epsilon^{-\frac{\lambda_1}{a}}
\eta^{\frac{\lambda_1}{a}} \lambda), \quad \eta \in (0,\delta)
\end{equation}
by Lemma~\ref{lm:matrix_de}. From~\eqref{f:critical_points_analytic_1}, we infer
that
\begin{equation*}     
 \tilde y = \lim_{\eta \downarrow 0} \varphi(\eta). 
\end{equation*}
In case C, Corollary~\ref{co:vanishing_flux} implies that $\tilde y = 0$. 
If $\nu$ was equal to $0$, it would then follow that $\varphi \equiv 0$ on
$(0,\delta)$. This is impossible in light of Lemma~\ref{lm:positive_densities}. 
As a result,
\begin{equation*}
 \varphi(\eta) = \nu \epsilon^{-\frac{\lambda_1}{a}} \eta^{\frac{\lambda_1}{a}}
\lambda + o(\eta^{\frac{\lambda_1}{a}})
\end{equation*}
as $\eta$ approaches $0$ from the right. This establishes part 2 of
Theorem~\ref{thm:critical_points_analytic} for case C and under the assumption
that $\tfrac{\lambda_1}{a}$ is not an integer. 

\medskip

In the resonance case, Proposition 11.6
in~\cite{Taylor_ode} implies that $Y^2 = 0$, that $Y \lambda = 0$ and that $Y
\tilde y$ is an eigenvector of $B_0$ corresponding to the eigenvalue
$\tfrac{\lambda_1}{a}$. Together with Lemma~\ref{lm:integer_eigenvalue}, this
yields 
\begin{align}      
\varphi(\eta) =& V(\eta)  \exp \Bigl( \ln \Bigl(
\frac{\eta}{\epsilon}
\Bigr) B_0 \Bigr)  \Bigl( \nu \lambda + \tilde y + \ln
\Bigl( \frac{\eta}{\epsilon} \Bigr)  Y (\nu \lambda + \tilde y)
\Bigr) \notag \\
=& V(\eta) 
\Bigl(\exp \Bigl( \ln \Bigl( \frac{\eta}{\epsilon} \Bigr) B_0 \Bigr)
(\nu \lambda + \tilde y) \label{f:critical_points_analytic_6} \\
& + \ln \Bigl( \frac{\eta}{\epsilon}
\Bigr)  \Bigl( Y \tilde y + \sum_{k=1}^{\infty} \frac{1}{k!} 
\Bigl( \ln \Bigl( \frac{\eta}{\epsilon} \Bigr) \Bigr)^k \Bigl(
\frac{\lambda_1}{a} \Bigr)^k Y \tilde y \Bigr) \Bigr).
\label{f:critical_points_analytic_3}
\end{align}
Using~\eqref{f:critical_points_analytic_4} and~\eqref{f:expansion_of_V}, the
term
in~\eqref{f:critical_points_analytic_6} and~\eqref{f:critical_points_analytic_3}
becomes 
\begin{equation}   \label{f:critical_points_analytic_5}
\Bigl( \mathbbm{1} + \sum_{k=1}^{\infty} \eta^k V_k \Bigr)  (\tilde y
+ \epsilon^{-\frac{\lambda_1}{a}} \eta^{\frac{\lambda_1}{a}}  (\nu \lambda
- \ln(\epsilon) Y \tilde y) + \epsilon^{-\frac{\lambda_1}{a}}
\eta^{\frac{\lambda_1}{a}} \ln(\eta) Y \tilde y).
\end{equation}
Let us first consider the situation where $\tfrac{\lambda_1}{a} > 1$. In case
C, $\tilde y = 0$ and we obtain 
\begin{equation*}
 \varphi(\eta) = \nu \epsilon^{-\frac{\lambda_1}{a}} \eta^{\frac{\lambda_1}{a}}
 \lambda + o(\eta^{\frac{\lambda_1}{a}})
\end{equation*}
as $\eta$ approaches $0$ from the right. 
Since $\nu \neq 0$ by Lemma~\ref{lm:positive_densities}, we have established
part 2 of Theorem~\ref{thm:critical_points_analytic} for case C under the
assumption that $\tfrac{\lambda_1}{a}$ is an integer larger than $1$. 

\medskip

Now, suppose that $\tfrac{\lambda_1}{a} = 1$. In case C,
Representation~\eqref{f:critical_points_analytic_5} of $\varphi(\eta)$
implies that 
\begin{equation*}
 \varphi(\eta) = \nu \epsilon^{-1} \eta  \lambda + o(\eta),
\end{equation*}
and part 3 of Theorem~\ref{thm:critical_points_analytic} follows for case
C.

In case B,~\eqref{f:critical_points_analytic_5} yields
\begin{equation*}
 \varphi(\eta) = \tilde y + \epsilon^{-1} \eta \ln(\eta) Y \tilde y + o(\eta
\ln(\eta)). 
\end{equation*}
Since $Y \tilde y$ is an eigenvector of $B_0$ with corresponding eigenvalue
$\tfrac{\lambda_1}{a}$, Lemma~\ref{lm:eigenvalues_B_0} implies that the
first component of $Y \tilde y$ is nonzero. This yields part 3 of
Theorem~\ref{thm:critical_points_analytic} for case B.
\epf

\bigskip

\section{Proof of Lemmas~\ref{lm:integral_equation_1}
and~\ref{lm:integral_equation_2}}   \label{sec:integral_equations}

For $t \geq 0$, as defined in \eqref{eq:invM} let $\mu \Pp^t$ denote
the distribution of $(X,A)_t$ starting from the initial distribution $\mu$. Since $\mu$ is invariant under
$(\Pp^t)_{t \geq 0}$, 
\begin{equation}       \label{f:invariant_measure}
\mu_i( ) = \int_{T^{(1)}}^{T^{(2)}} \pi(t)
 \mu \Pp^t( \times
\{i\}) \,dt
\end{equation} 
for any $T^{(1)} < T^{(2)}$ in $[0,\infty]$ and for any probability
density $\pi(t)$ on $(T^{(1)},T^{(2)})$. 

We will expand the expression on the right with respect to the sequences
of driving
vector fields and will ultimately see how $\rho_i$
gets transformed through the action of the Markov semigroup and through
time-averaging.

The following formula is the key to Lemmas~\ref{lm:integral_equation_1}
and~\ref{lm:integral_equation_2}.

\begin{lemma}  \label{lm:probability_densities}
Let $E \subset \R$ be a Borel set and let $i \in S$. For any
$T^{(1)} < T^{(2)}$ in $[0,\infty]$ and for any probability density $\pi(t)$ on
$(T^{(1)},T^{(2)})$, 
\begin{multline*}
\mu_i(E) 
= \int_E \Bigl( \int_{T^{(1)}}^{T^{(2)}} \pi(t)  e^{-\lambda_i t}
 \,\Phi_i^t
\# \rho_i(\eta) \,dt \\
 + \sum_{j \neq i} \lambda_{j,i}  \Bigl( \int_0^{T^{(1)}} e^{-\lambda_i
t} 
\,\Phi_i^t \# \rho_j(\eta) \,dt + \int_{T^{(1)}}^{T^{(2)}} c(t) 
e^{-\lambda_i
t} \, \Phi_i^t \# \rho_j(\eta) \,dt \Bigr) \Bigr) \,d\eta,
\end{multline*}
where $c(t):=\int_0^{T^{(2)}-t} \pi(s+t) \,ds$. 
\end{lemma}

Given this representation for $\mu_i$, we first show
Lemma~\ref{lm:integral_equation_1} and then
Lemma~\ref{lm:integral_equation_2}. Finally, we prove the representation
itself. 

\subsection{Proof of Lemma~\ref{lm:integral_equation_1}}

When we set $T^{(1)}=0$, $T^{(2)} = T$ and $\pi(t) = \tfrac{1}{T}$, the
identity in Lemma~\ref{lm:probability_densities} becomes
\begin{multline*}
\mu_i(E) = \int_E \Bigl(\frac{1}{T} \int_0^T e^{-\lambda_i t}
\, \Phi_i^t
\# \rho_i(\eta) \, dt \\
 + \sum_{j \neq i} \lambda_{j,i}  \frac{1}{T} \int_0^T (T-t) 
e^{-\lambda_i t} 
\,\Phi_i^t \# \rho_j(\eta) \,dt \Bigr) \, d\eta.
\end{multline*}
This implies Lemma~\ref{lm:integral_equation_1}. 

\subsection{Proof of Lemma~\ref{lm:integral_equation_2}}

When we set $T^{(1)}=T$ for some time $T > 0$, $T^{(2)} = \infty$ and $\pi(t) =
e^{T-t}$, the identity in Lemma~\ref{lm:probability_densities} becomes 
\begin{align*}
\mu_i(E) =& \int_E \Bigl( \int_{T}^{\infty} e^{T-t}  e^{-\lambda_i t}
 \,\Phi_i^t
\# \rho_i(\eta) \,dt \\
& + \sum_{j \neq i} \lambda_{j,i}  \Bigl( \int_0^{T} e^{-\lambda_i
t} 
\,\Phi_i^t \# \rho_j(\eta) \,dt + \int_{T}^{\infty} e^{T-t} 
e^{-\lambda_i
t} \, \Phi_i^t \# \rho_j(\eta) \,dt \Bigr) \Bigr) \,d\eta.
\end{align*}
Since $\mu$ is a probability measure, 
\begin{align*}
\int_E \int_T^{\infty} e^{T-t}  e^{-\lambda_i t} \, \Phi_i^t \#
\rho_i(\eta) \,dt \,d\eta =& \int_T^{\infty} e^{T-t}  e^{-\lambda_i t}
 \mu_i((\Phi_i^t)^{-1}(E)) \,dt \\
\leq & e^{-\lambda_i T}  \int_T^{\infty} e^{T-t} \,dt = e^{-\lambda_i T}, 
\end{align*}
where one should observe that the set $(\Phi_i^t)^{-1}(E)$ is well-defined even
if $\Phi_i^{-t}(\eta)$ is undefined for some $\eta \in E$. 
Similarly, 
\begin{equation*}
\sum_{j \neq i} \lambda_{j,i}  \int_E \int_T^{\infty} e^{T-t} 
e^{-\lambda_i t}\,  \Phi_i^t \# \rho_j(\eta) \,dt \,d\eta \leq  \sum_{j \neq
i} \lambda_{j,i}  e^{-\lambda_i T}.
\end{equation*}
Letting $T$ go to infinity, we obtain that
\begin{equation*}
 \mu_i(E) = \int_E \sum_{j \neq i} \lambda_{j,i} 
\int_0^{\infty} e^{-\lambda_i t} \, \Phi_i^t \# \rho_j(\eta) \,dt \,d\eta,
\end{equation*}
and Lemma~\ref{lm:integral_equation_2} follows. 

\subsection{Proof of Lemma~\ref{lm:probability_densities}}

Fix an $i \in S$. We introduce some notation. For any $t>0$ and
for any index sequence $\ibf$ with terminal
index $i$, let $C^t_{\ibf}$ denote
the
event that the driving vector fields up to time $t$ appear in the order given
by $\ibf$. For any index sequence $\ibf =(i_1, \ldots, i_{m-1},i)$ of length
$m \geq 2$, let
$\Pp_{\ibf}$ be the probability that the first $m$ driving vector fields appear
in the order given by $\ibf$, conditioned on $u_{i_1}$ being the first driving
vector field. 
For $T > 0$ and $m \in \N$, we define the simplex $\Delta_{T,m}$ as
the interior of the
convex hull of the vectors $0, T e_1, \ldots, T e_m$ in $\R^m$. For any vector
$v$
with $m$ components, no matter whether $v$ is a vector of indices, switching
times or switching rates, let $v^{(m-1)}$ denote the projection of $v$ onto its
first $(m-1)$ coordinates. Moreover, let $\| v \|_1$ be the sum of the
coordinates of $v$ and let $\langle \cdot, \cdot \rangle$ denote the Euclidean
inner product on the space that fits the context (usually $\R^{m-1}$ or
$\R^m$).

\begin{lemma}   \label{lm:one_switch}
For any $T^{(1)} < T^{(2)}$ in $[0, \infty]$ and for any function $\pi(t)$
that is nonnegative and integrable on $(T^{(1)}, T^{(2)})$,  
\begin{multline*}
 \int_{T^{(1)}}^{T^{(2)}} \pi(t)  \Pp_{\xi,i}(C^t_{(i)})  \int_{\R}
\Pp_{\xi,i}(X_t \in E \vert C^t_{(i)}) \,\mu_i(d\xi) \,dt \\
= \int_E \int_{T^{(1)}}^{T^{(2)}} \pi(t)  e^{-\lambda_i t} \, \Phi_i^t
\# \rho_i(\eta) \,dt \,d\eta.
\end{multline*}
\end{lemma}

\bpf This is immediate.
\epf

\bigskip

\begin{lemma} \label{lm:identity}
For any index sequence $\ibf=(i_1, \ldots, i_{m-1},i)$ of length $m
\geq 2$, for any $T^{(1)} < T^{(2)}$ in $[0,\infty]$ and for any function
$\pi(t)$ that is nonnegative and integrable on $(T^{(1)},T^{(2)})$, 
\begin{multline*}
 \int_{T^{(1)}}^{T^{(2)}} \pi(t) 
\Pp_{\xi,i_1}(C^t_{\ibf})  \int_{\R}
\Pp_{\xi,i_1}(X_t \in E \vert C^t_{\ibf}) \,\mu_{i_1}(d\xi) \,dt
\\
= \Pp_{\ibf} \prod_{l=1}^{m-1} \lambda_{i_l} \int_{\Delta_{T^{(2)},m} \setminus
\Delta_{T^{(1)},m}}
\pi(\| \sbf \|_1)  e^{- \langle \lambdabf^{(m-1)}, \sbf^{(m-1)} \rangle}
 e^{-\lambda_i s_m}  \mu_{i_1}((\Phi_{\ibf}^{\sbf})^{-1}(E)) \,d\sbf,
\end{multline*}
where $\lambdabf^{(m-1)} = (\lambda_{i_1}, \ldots, \lambda_{i_{m-1}})^T$. 
\end{lemma}

\bpf Fix an index sequence $\ibf = (i_1, \ldots, i_{m-1}, i)$ of length $m$,
$T^{(1)} < T^{(2)} \in [0,\infty]$ and a nonnegative integrable function $\pi$
on $(T^{(1)}, T^{(2)})$. Let $T_1, \ldots, T_m$ be independent, exponentially
distributed random variables such
that $T_l$ has parameter $\lambda_{i_l}$ for $1 \leq l \leq m-1$ and
$T_m$ has parameter $\lambda_i$. For any $t \geq T^{(1)}$,  
\begin{multline}
  \int_{\R} \Pp_{\xi,i_1}(X_t \in E \vert C^t_{\ibf}) \,\mu_{i_1}(d\xi)
 \\
= \frac{1}{\Pp(R^t_{\ibf})}  \int_{\R} \Pp(\Phi_{\ibf}^{(T_1, \ldots,
T_{m-1}, t-\sum_{l=1}^{m-1} T_l)}(\xi) \in E, R^t_{\ibf}) \
\mu_{i_1}(d\xi),    \label{f:identity} 
\end{multline}
where
\begin{equation*}
 R^t_{\ibf} := \Bigl\{ \sum_{l=1}^{m-1} T_l < t \leq \sum_{l=1}^m T_l \Bigr
\}.
\end{equation*}
As a notational shorthand, we introduce the functions
\begin{equation*}
 f^{\xi}_{t,\ibf}\colon \,\R^{m-1} \to \R, \,(s_1, \ldots, s_{m-1}) \mapsto
\Phi_{\ibf}^{(s_1, \ldots, s_{m-1}, t - \sum_{l=1}^{m-1} s_l)}(\xi).
\end{equation*}
Then,
\begin{multline*}
 \Pp(\Phi_{\ibf}^{(T_1, \ldots, T_{m-1}, t-\sum_{l=1}^{m-1} T_l)}(\xi) \in E;
R^t_{\ibf}) \\
= \int_{\Delta_{t,m-1}} \mathbbm{1}_{\{ f^{\xi}_{t,
\ibf}(\sbf) \in E \}}(\sbf)  \prod_{l=1}^{m-1} \lambda_{i_l} 
e^{-\lambda_{i_l} s_l}  e^{-\lambda_i (t -\| \sbf \|_1)} \,d\sbf,
\end{multline*}
which implies that~\eqref{f:identity} can be written as 
\begin{equation}     \label{f:identity_2}
\frac{1}{\Pp(R^t_{\ibf})}  \int_{\R} \int_{\Delta_{t,m-1}}
\mathbbm{1}_{\{f^{\xi}_{t,\ibf}(\sbf) \in
E\}}(\sbf)  \prod_{l=1}^{m-1} \lambda_{i_l} 
e^{-\lambda_{i_l} s_l}  e^{-\lambda_i (t - \| \sbf
\|_1)} \,d\sbf \,\mu_{i_1}(d\xi).
\end{equation}
Interchanging the order of integration,the righthand side of ~\eqref{f:identity_2} becomes
\begin{equation*}
 \frac{1}{\Pp(R^t_{\ibf})}  \int_{\Delta_{t,m-1}}
\prod_{l=1}^{m-1} \lambda_{i_l}  e^{-\lambda_{i_l} s_l-\lambda_i (t - \| \sbf \|_1)} 
\mu_{i_1}((f^{\cdot}_{t,\ibf}(\sbf))^{-1}(E)) \,d\sbf.
\end{equation*}
We have thus shown that
\begin{multline*}
 \int_{T^{(1)}}^{T^{(2)}} \pi(t) 
\Pp_{\xi,i_1}(C^t_{\ibf})  \int_{\R}
\Pp_{\xi,i_1}(X_t \in E \vert C^t_{\ibf}) \,\mu_{i_1}(d\xi) \, dt
= \int_{T^{(1)}}^{T^{(2)}} \pi(t) 
\frac{\Pp_{\xi,i_1}(C^t_{\ibf})}{\Pp(R^t_{\ibf})} \\ \cdot
 \int_{\Delta_{t,m-1}} \prod_{l=1}^{m-1}
\lambda_{i_l}  e^{-\lambda_{i_l} s_l-\lambda_i (t - \| \sbf \|_1)} 
\mu_{i_1}((f^{\cdot}_{t,\ibf}(\sbf))^{-1}(E)) \,
d\sbf \, dt.
\end{multline*}
The term $\tfrac{\Pp_{\xi,i_1}(C^t_{\ibf})}{\Pp(R^t_{\ibf})}$ gives
the probability that the first $m$ driving vector fields appear according to
index sequence $\ibf$, conditioned on $u_{i_1}$ being the first driving vector
field. It is
therefore equal to $\Pp_{\ibf}$.   
Interchanging the order of integration and
substituting $s_m =
t-\| \sbf \|_1$, Lemma~\ref{lm:identity} follows. 
\epf

\medskip

\begin{lemma} \label{lm:time_after_last_switch}
For any index sequence $\ibf=(i_1, \ldots, i_{m-1},i)$ of length $m
\geq 2$, for any $T^{(1)} < T^{(2)}$ in $[0,\infty]$ and for any
function $\pi(t)$ that is nonnegative and integrable on $(T^{(1)},T^{(2)})$,
\begin{multline*}
\int_{T^{(1)}}^{T^{(2)}} \pi(t)  \Pp_{\xi,i_1}(C^t_{\ibf}) 
\int_{\R}
\Pp_{\xi,i_1}(X_t \in E \vert C^t_{\ibf}) \,\mu_{i_1}(d\xi) \,dt 
= \\\int_{\Delta_{T^{(2)},2} \setminus \Delta_{T^{(1)},2}} \lambda_{i_{m-1},i}
e^{-\lambda_i t}  \pi(s+t)  \Pp_{\xi,i_1}(C^s_{\ibf^{(m-1)}})  \\
\cdot \int_{\R} \Pp_{\xi,i_1}(X_s \in (\Phi_i^t)^{-1}(E) \vert
C^s_{\ibf^{(m-1)}}) \,\mu_{i_1}(d\xi) \,d(s,t).
\end{multline*}
\end{lemma}

\bpf For notational compactness momentarily we introduce the notation $\tilde
\Delta_{i}(m,t):=\Delta_{T^{(i)} - t, m-1}$. Then by Tonelli's theorem,
the term to the right of the equality sign in Lemma~\ref{lm:identity}
can be written as
\begin{multline*}
\int_0^{T^{(1)}} \Pp_{\ibf}  \prod_{l=1}^{m-1} \lambda_{i_l} 
\int_{(\tilde \Delta_{2}\setminus \tilde \Delta_{1})(m,t)}
\pi( \| \sbf \|_1 + t)  e^{- \langle \lambdabf^{(m-1)}, \sbf
\rangle-\lambda_i t}  
  \mu_{i_1}((\Phi_{\ibf}^{(\sbf,t)})^{-1}(E)) \,
d\sbf \,dt  \\
 + \int_{T^{(1)}}^{T^{(2)}} \Pp_{\ibf}  \prod_{l=1}^{m-1}
\lambda_{i_l} 
\int_{\tilde \Delta_2(m,t)} \pi(
\| \sbf \|_1 + t)  e^{- \langle \lambdabf^{(m-1)}, \sbf
\rangle}  e^{-\lambda_i t} 
  \\ \cdot \mu_{i_1}((\Phi_{\ibf}^{(\sbf,t)})^{-1}(E)) \,d\sbf \,dt
 \\
= \int_0^{T^{(1)}} \lambda_{i_{m-1},i} 
e^{-\lambda_i t}  \Pp_{\ibf^{(m-1)}}  \prod_{l=1}^{m-2}
\lambda_{i_l}  \int_{(\tilde \Delta_{2}\setminus \tilde \Delta_{1})(m,t)} 
\pi_t( \| \sbf \|_1)  
  e^{- \langle
\lambda^{(m-1)}, \sbf \rangle} \\ \cdot
\mu_{i_1}((\Phi_{\ibf^{(m-1)}}^{\sbf})^{-1}((\Phi_i^t)^{-1}(E))) \,d\sbf \,dt
\\
 + \int_{T^{(1)}}^{T^{(2)}} \lambda_{i_{m-1},i}  e^{-\lambda_i t}
\Pp_{\ibf^{(m-1)}}  \prod_{l=1}^{m-2} \lambda_{i_l}
\int_{\tilde \Delta_2(m,t)} \pi_t( \| \sbf \|_1)   
 e^{- \langle \lambdabf^{(m-1)}, \sbf \rangle} \\
\cdot \mu_{i_1}((\Phi_{\ibf^{(m-1)}}^{\sbf})^{-1}((\Phi_i^t)^{-1}(E))) \,d\sbf
\,dt,   
\end{multline*}
where the function $\pi_t(s) := \pi(s+t)$ is nonnegative and integrable on
$(T^{(1)} -t,T^{(2)}-t)$ if $t < T^{(1)}$, and is nonnegative and
integrable on
$(0,T^{(2)}-t)$ if $t > T^{(1)}$. 

By another application of Lemma~\ref{lm:identity} (if $m >2$) or of
Lemma~\ref{lm:one_switch} (if $m=2$),
the previous term becomes
\begin{multline*}
\int_0^{T^{(1)}} \lambda_{i_{m-1},i}  e^{-\lambda_i t} 
\int_{T^{(1)}-t}^{T^{(2)}-t} \pi_t(s)  \Pp_{\xi,i_1}(C^s_{\ibf^{(m-1)}}) \\\cdot
 \int_{\R} \Pp_{\xi,i_1}(X_s \in (\Phi_i^t)^{-1}(E) \vert
C^s_{\ibf^{(m-1)}}) \,\mu_{i_1}(d\xi) \,ds \,dt \\
 + \int_{T^{(1)}}^{T^{(2)}} \lambda_{i_{m-1},i}  e^{-\lambda_i t}
 \int_0^{T^{(2)}-t} \pi_t(s)
 \Pp_{\xi,i_1}(C^s_{\ibf^{(m-1)}}) \\\cdot
  \int_{\R} \Pp_{\xi,i_1}(X_s \in (\Phi_i^t)^{-1}(E) \vert
C^s_{\ibf^{(m-1)}}) \,\mu_{i_1}(d\xi) \,ds \,dt \\
= \int_{\Delta_{T^{(2)},2} \setminus \Delta_{T^{(1)},2}} \lambda_{i_{m-1},i}
e^{-\lambda_i t}  \pi(s+t)  \Pp_{\xi,i_1}(C^s_{\ibf^{(m-1)}})
\\\cdot
  \int_{\R} \Pp_{\xi,i_1}(X_s \in (\Phi_i^t)^{-1}(E) \vert
C^s_{\ibf^{(m-1)}}) \,\mu_{i_1}(d\xi) \,d(s,t).
\end{multline*}
\epf

\bpf [Proof of
Lemma~\ref{lm:probability_densities}] Fix a Borel set $E$, $T^{(1)} < T^{(2)}
\in [0,\infty]$ and a probability
density $\pi$ on $(T^{(1)},T^{(2)})$. Expanding the
term to the right of the equality sign in~\eqref{f:invariant_measure} by
conditioning on the sequences of driving vector fields, we obtain
\begin{align}
\mu_i(E) =& \int_{T^{(1)}}^{T^{(2)}} \pi(t)  \Pp_{\xi,i}(C^t_{(i)})
 \int_{\R}
\Pp_{\xi,i}(X_t \in E \vert C^t_{(i)}) \,\mu_i(d\xi) \,dt \notag \\
& + \sum_{\ibf: \lvert \ibf \rvert \geq 2}  \int_{T^{(1)}}^{T^{(2)}}
\pi(t) 
\Pp_{\xi,i_1}(C^t_{\ibf})  \int_{\R} \Pp_{\xi,i_1}(X_t \in E \vert
C^t_{\ibf}) \,\mu_{i_1}(d\xi) \,dt,   \label{f:probability_densities}
\end{align}
where $\sum_{\ibf: \lvert \ibf \rvert \geq 2}$ extends over all index
sequences $\ibf = (i_1, \ldots, i_{m-1}, i)$ with terminal index $i$ and length
$m \geq 2$. 

By Lemma~\ref{lm:one_switch}, it is enough to show that the term
in~\eqref{f:probability_densities} equals 
\begin{equation*}
\sum_{j \neq i} \lambda_{j,i}  \int_E \Bigl( \int_0^{T^{(1)}}
e^{-\lambda_i
t} \,
\Phi_i^t \# \rho_j(\eta) \,dt + \int_{T^{(1)}}^{T^{(2)}} c(t) 
e^{-\lambda_i
t}  \Phi_i^t \# \rho_j(\eta) \,dt \Bigr) \
d\eta,
\end{equation*}
where $c(t)$ is defined as in Lemma~\ref{lm:probability_densities}.
For any $m \geq 2$, let $\sum_{\ibf: \lvert \ibf \rvert = m}$ be the sum over
all index sequences of length $m$ with terminal index $i$. For any $j \in S$,
let $\sum_{\ibf}^{(j)}$ be the sum over all index
sequences $\ibf$ with terminal index $j$. By
Lemma~\ref{lm:time_after_last_switch}, the term
in~\eqref{f:probability_densities} can be written as  
\begin{multline}\label{f:probability_densities-2}
 \sum_{m=2}^{\infty} \sum_{\ibf: \lvert \ibf \rvert = m}
\int_{\Delta_{T^{(2)},2} \setminus
\Delta_{T^{(1)},2}}
\lambda_{i_{m-1},i}  e^{-\lambda_i t}  \pi(s+t)
\Pp_{\xi,i_1}(C^s_{\ibf^{(m-1)}})\\ \cdot
 \int_{\R} \Pp_{\xi,i_1}(X_s \in (\Phi_i^t)^{-1}(E) \vert
C^s_{\ibf^{(m-1)}}) \,\mu_{i_1}(d\xi) \,d(s,t)  \\
= \sum_{j \neq i} \int_{\Delta_{T^{(2)},2} \setminus \Delta_{T^{(1)},2}}
\lambda_{j,i} 
e^{-\lambda_i t}  \pi(s+t) \\\cdot
  \sum_{\ibf}^{(j)} \Pp_{\xi,i_1}(C^s_{\ibf})
 \int_{\R} \Pp_{\xi,i_1}(X_s \in (\Phi_i^t)^{-1}(E) \vert C^s_{\ibf}) \
\mu_{i_1}(d\xi) \,d(s,t).  
\end{multline}
Moreover, for any fixed $s$,  
\begin{equation*}
\sum_{\ibf}^{(j)} \Pp_{\xi,i_1}(C^s_{\ibf}) 
\int_{\R} \Pp_{\xi,i_1}(X_s \in (\Phi_i^t)^{-1}(E) \vert C^s_{\ibf}) \
\mu_{i_1}(d\xi) 
= \mu \Pp^s((\Phi_i^t)^{-1}(E) \times \{j\}). 
\end{equation*}
Since $\mu$ is invariant, $\mu
\Pp^s((\Phi_i^t)^{-1}(E) \times \{j\})$ equals $\mu_j((\Phi_i^t)^{-1}(E))$
and is thus independent of $s$.

As a result, the right side of~\eqref{f:probability_densities-2} equals
\begin{multline*}
\sum_{j \neq i} \int_{\Delta_{T^{(2)},2} \setminus \Delta_{T^{(1)},2}}
\lambda_{j,i} 
e^{-\lambda_i t}  \pi(s+t)  \mu_j((\Phi_i^t)^{-1}(E)) \
d(s,t) \\
= \sum_{j \neq i} \int_0^{T^{(1)}} \lambda_{j,i}  e^{-\lambda_i t} 
\mu_j((\Phi_i^t)^{-1}(E))  \int_{T^{(1)}-t}^{T^{(2)}-t}
\pi(s+t) \
ds \,dt \\
 + \sum_{j \neq i} \int_{T^{(1)}}^{T^{(2)}} \lambda_{j,i} 
e^{-\lambda_i t} 
 \mu_j((\Phi_i^t)^{-1}(E))  \int_0^{T^{(2)} -t}
\pi(s+t)
\,ds \,dt,
\end{multline*}
and Lemma~\ref{lm:probability_densities} follows. 
\epf 

\bigskip

\appendix 

\section{How Equation~\eqref{f:representation_infinite_time} relates to the
Fokker--Planck equations}   \label{sec:Fokker_Planck}

Equation~\eqref{f:representation_infinite_time} in the proof of
Lemma~\ref{lm:representation_infinite_time} can also be derived from the
Fokker--Planck equations, 
but in order to do this, one needs to assume that the invariant densities are
$\Cs^1$.

It is an immediate consequence of~\eqref{f:Fokker_Planck} that
\begin{equation*}
\bar \rho(\zeta) = (\lambda_1 + u_1'(\zeta))  \rho_1(\zeta) + u_1(\zeta)
 \rho_1'(\zeta),
\end{equation*}
see~\cite{Balazs}. Hence, the term to the right of the equality sign
in~\eqref{f:representation_infinite_time} equals
\begin{align}
& -\frac{1}{u_1(\eta)}  \int_{\eta}^{\vartheta} (\lambda_1 +
u_1'(\zeta))  \rho_1(\zeta)  \exp \Bigl( \lambda_1 
\int_{\eta}^{\zeta} \frac{dx}{u_1(x)} \Bigr) \,d\zeta \label{f:FP1} \\
& -\frac{1}{u_1(\eta)}  \int_{\eta}^{\vartheta} \rho_1'(\zeta) 
u_1(\zeta)  \exp \Bigl( \lambda_1  \int_{\eta}^{\zeta}
\frac{dx}{u_1(x)} \Bigr) \,d\zeta. \label{f:FP2} 
\end{align}
As 
\begin{equation*}
 \lim_{\zeta \uparrow \vartheta} \Bigl( \rho_1(\zeta)  u_1(\zeta) 
\exp \Bigl( \lambda_1  \int_{\eta}^{\zeta} \frac{dx}{u_1(x)} \Bigr)
\Bigr) = 0
\end{equation*}
if $u_1$ is smooth and forward-complete, integration by parts implies that the
term in~\eqref{f:FP2} equals
\begin{equation}  \label{f:FP3}
\rho_1(\eta) + \frac{1}{u_1(\eta)}  \int_{\eta}^{\vartheta}
(\lambda_1 + u_1'(\zeta))  \rho_1(\zeta) 
\exp \Bigl( \lambda_1 
\int_{\eta}^{\zeta} \frac{dx}{u_1(x)} \Bigr) \,d\zeta. 
\end{equation}
Since the second term in~\eqref{f:FP3} cancels with the term in~\eqref{f:FP1},
we obtain~\eqref{f:representation_infinite_time}. 
\bibliographystyle{alpha}

\bibliography{regularity}

\end{document}